\documentclass[preprint,12pt]{elsarticle}

\usepackage[letterpaper,top=2cm,bottom=2cm,left=3cm,right=3cm,marginparwidth=1.75cm]{geometry}

\usepackage{amsfonts,amsmath,amssymb}
\usepackage{tikz}
\usepackage[colorlinks=true, citecolor=black, linkcolor=black, urlcolor=black]{hyperref}

\newtheorem{thm}{Theorem}[section]
\newtheorem{lem}[thm]{Lemma}

\newtheorem{prb}[thm]{Problem}
\newtheorem{prop}[thm]{Proposition}
\newtheorem{df}{Definition}
\newtheorem{rem}{Remark}

\date{}
\title{
The Lawrence-Krammer  representation is a quantization of the symmetric square of
the Burau representation
\footnote{To all fearless Ukrainians defending not only their country, but the whole civilization 
against putin's 
\href{https://en.wikipedia.org/wiki/Ruscism}{rashism}
}
}
\author{A.V.~Kosyak}
\address{Institute of Mathematics, Ukrainian National Academy of Sciences,
3 Tereshchenkivs'ka Str., Kyiv, 01601, Ukraine}
\address{Royal Institution,
London Institute for Mathematical Sciences, 21 Albemarle St, London W1S 4BS, UK}
\journal{
}
\begin{document}

\maketitle

\begin{abstract}
We show that the Lawrence--Krammer  representation 
%can be obtained as 
is the quantization of the symmetric square  of  the Burau representation. This connection allows us to construct new representations of braid groups.
\end{abstract}

\iffalse
{\it keywords}:
braid group,  Lawrence-Krammer  representation, Burau representation,
symmetric tensor power, quantization,  quantum group, $q-$Pascal triangle,
knot invariant

%\begin{keyword}
%Quantum groups (quantized function algebras) and their representations
%\end{keyword}
\fi

\tableofcontents
\section{Introduction}
In 1936  W.~Burau  introduced \cite{Bur36} the family of representations of the braid group $B_n$ as the
{\it deformation}
%({\it not quantization})
of the {\it  standard representation} of the symmetric group $S_n$ (see (\ref{Bur_n})).
 %For a long time it was believed to be faithful.
The Burau representation for $n = 2, 3$ has been known to be faithful for some time.
Indeed,  in 1969 Magnus and Peluso \cite{MagPel69} showed that the Burau and {\it Gasner representations} are faithful for $n=2,3$. See also J.~Birman \cite[Theorem 3.15]{Bir74}.
Moody \cite{Mod91} in 1991 showed that the Burau representation is not
faithful for $n \geq 9$; this result was improved to $n \geq 6$ by Long and Paton in 1992, \cite{LonPat93}. The non-faithfulness for $n = 5$ was shown by Bigelow in 1999, \cite{Big99}. The faithfulness of the Burau representation for $n = 4$ is an {\it open problem}.

In 1990 R.~Lawrence \cite{Law90} have constructed two-parameter family of representation of the braid group $B_n$, using the homological methods. In 2000 D.~Krammer  \cite{Kra00}  using completely algebraic methods have constructed the same representations
%equivalent  family of two parameter representations of the Braid group $B_n$
and showed that this representation is faithful for $B_4$. One year later, in 2001, S.~Bigelow showed that the Lawrence-Krammer representation is faithful for all $B_n$. Next year, D.~Krammer \cite{Kra02} proved the same result by a different method.

The Burau representation appears as a summand of the {\it Jones representation} \cite{Jon87}, and for $n = 4$, the faithfulness of the Burau representation is equivalent to that of the Jones representation, which on the other hand is related to the question of whether or not the Jones polynomial is an {\it unknot detector} \cite{Big02.11}.

In \cite{JacKel11} Jackson and Kerler established an isomorphism between Lawrence-Krammer representation
%gave an explicit construction and proof of an isomorphism between the
%faithful representation $H_{n,2}$ of $B_n$ considered by Bigelow and Krammer
and the submodule
of the $R$-matrix representations on $V^{\otimes n}$ for the {\it generic Verma module} $V$ of the quantum
group $U_q (\mathfrak{sl}_2)$.
%
%\begin{rem} \label{r.Kram-other}

From  \cite{JacKel11}: ``In \cite{Zin01} Zinno manages to find a different identification of the Lawrence-Krammer representation with a
quantum algebraic object, namely the quotient of the Birman-Wenzl-Murakami algebra similarly
defined over ${\mathbb Z}[q^\pm,t^\pm]$. This representation can, by \cite{Wenz90}, be understood as the one arising
from the quantum orthogonal groups $U_q(\mathfrak{so}(k + 1))$ acting on the $n$-fold tensor product of the
fundamental representation.'' See \cite[n 4.6, p.71]{BirBre93} for more details about connections between the Lawrence-Krammer and other representations of braid groups.

\section{Connection between the Lawrence--Krammer and the Burau representations}
\subsection{The Burau representation}
 {\it  The Burau representation}
$\rho:B_{n}\to {\rm GL}_n({\mathbb Z}[t,t^{-1}])$  is defined
%(see c.f. \cite{Jon87} )
for  a non-zero complex number $t$ by
\begin{equation}
\label{Bur_n}
\sigma_i\mapsto
I_{i-1}\oplus\left(\begin{smallmatrix}
1-t&t\\
1&0\\
\end{smallmatrix}\right)\oplus I_{n-i-1}
\end{equation}
where $1-t$ is the $(i,i)$ entry.
For $t=1$ this is the {\it standard representation} of the symmetric group $S_n$
interchanging by matrix
$
\left(\begin{smallmatrix}
0&1\\
1&0\\
\end{smallmatrix}\right)
$
two neighbors basic elements $e_i$ and $e_{i+1}$  in the space ${\mathbb C}^n$.
The vector $e=e_1+\dots+e_n$ is invariant, therefore
the representation $\rho$ splits into 1-dimensional and an
$n\!-\!1-$dimensional irreducible representations, known as the {\it
reduced Burau representation} $\rho_n^{(t)}:B_{n}\to {\rm
GL}_{n-1}({\mathbb Z}[t,t^{-1}])$ %%
\begin{equation}
\label{r-Bur_n}
\sigma_1\mapsto
% b_1=
 \left(\begin{smallmatrix}
-t&0\\
-1&1\\
\end{smallmatrix}\right)\oplus I_{n-3},\quad
\sigma_{n-1}\mapsto
I_{n-3}\oplus \left(\begin{smallmatrix}
1&-t\\
0&-t\\
\end{smallmatrix}\right),
\end{equation}
\begin{equation}
\label{r-Bur_n(b)}
\sigma_i\mapsto
%b_i=
I_{i-2}\oplus\left(\begin{smallmatrix}
1&-t&0\\
0&-t&0\\
0&-1&1\\
\end{smallmatrix}\right)\oplus I_{n-i-2},\,\,2\leq i\leq n-2.
\end{equation}
We use the following equivalent form of  the
reduced Burau representation $\rho_n^{(t)}:B_{n}\mapsto {\rm
GL}_{n-1}({\mathbb Z}[t,t^{-1}])$, compare with (\ref{r-Bur_n}) and (\ref{r-Bur_n(b)}):
\begin{equation}
\label{Bur_n-conj}
 \sigma_1\mapsto
 %b_1=
 \left(\begin{smallmatrix}
-t&t\\
0&1\\
\end{smallmatrix}\right)\oplus I_{n-3},\quad
\sigma_{n-1}\mapsto
%b_{n-1}=
I_{n-3}\oplus \left(\begin{smallmatrix}
1&0\\
1&-t\\
\end{smallmatrix}\right),
\end{equation}
\begin{equation}
\label{Bur_n(b)-conj} \sigma_i\mapsto
I_{i-2}\oplus\left(\begin{smallmatrix}
1&0&0\\
1&-t&t\\
0&0&1\\
\end{smallmatrix}\right)\oplus I_{n-i-2},\,\,2\leq i\leq n-2.
\end{equation}
The equivalence  follows from the following:
\begin{equation*}
 J^{-1}
\left(\begin{smallmatrix}
1&-t&0\\
0&-t&0\\
0&-1&1\\
\end{smallmatrix}\right)
J=
\left(\begin{smallmatrix}
1&1&0\\
0&-t&0\\
0&t&1\\
\end{smallmatrix}\right)\stackrel{(\cdot)^t}{\mapsto}
\left(\begin{smallmatrix}
1&0&0\\
1&-t&t\\
0&0&1\\
\end{smallmatrix}\right)
,\quad\text{where}\quad 
J=\left(\begin{smallmatrix}
0&0&1\\
0&-1&0\\
1&0&0\\
\end{smallmatrix}\right).
\end{equation*}
\subsection{The Burau representation of  $B_\infty$}
\label{ch.B(inf)}
Define the group $B_\infty$ as follows
\begin{equation}
\label{B(inf)}
B_\infty=\langle(\sigma_i)_{i\,\in\,\mathbb Z}
\mid\sigma_i\sigma_{i+1}\sigma_i=\sigma_{i+1}\sigma_i\sigma_{i+1},\,\,i\in{\mathbb Z},\,\,\,\,
 \sigma_i\sigma_j=\sigma_i\sigma_j,\quad\mid i-j\mid\,\geq
2\rangle.
\end{equation}
Consider a complex Hilbert space
\begin{equation}
\label{l_2(Z)}
H=l_2(\mathbb Z)=\{x=(x_k)_{k\,\in\,\mathbb Z}\mid \Vert x\Vert^2=\sum_{k\in\mathbb Z}\vert x_k\vert^2<\infty\}
\end{equation}
with its standard orthonormal basis $e=(e_k)_{k\,\in\,\mathbb Z}$ defined by $e_k=(\delta_{k,n})_{n\,\in\,\mathbb Z}$.  Define the reduced Burau representation of the group
$B_\infty$ in the space $H$ as follows
\begin{equation}
\label{B(inf)-rep-bur}
\sigma_i\mapsto
b_{i,\infty}=I_{\infty}\oplus\left(\begin{smallmatrix}
1&0&0\\
1&-t&t\\
0&0&1\\
\end{smallmatrix}\right)\oplus I_{\infty},\quad i\in \mathbb Z,
\end{equation}
where $-t$ is the $(i,i)$ entry. It is clear that $b_{i+1,\infty}$ is obtained from $b_{i,\infty}$ by shifting  the basis $e: e_n\to e_{n+1}$. More precisely, if we denote by $J$ the unitry operator on $l_2(\mathbb Z)$ defined as follows: $Je_n= e_{n+1},\,\,n\in \mathbb Z$ we get
\begin{equation}
\label{J(inf)}
Jb_{i,\infty}J^{-1}=b_{i+1,\infty}.
\end{equation}
Next, denote by $P_n$ the orthogonal projector of $H$ on the
subspace $H_n$ generated by the first $n$ vectors, i.e., $H_n=\langle e_k\mid 1\leq k\leq n\rangle$. Then we get
\begin{equation}
\label{B(inf)-red-Bur}
P_nb_{k,\infty}P_n= \rho_{n+1}^{(t)}(\sigma_k),\quad 1\leq k\leq n.
\end{equation}
Denote by $J_n$ the unitary operator on $H_n$ defined as follows
\begin{equation}
\label{J_n}
J_ne_k=e_{k+1},\quad 1\leq k \leq n-1,\quad J_ne_n=e_{1},
\end{equation}
and by $i_n$ the natural embedding of ${\rm GL}(n-1,\mathbb C)$ into ${\rm GL}(n,\mathbb C)$ defined as follows:
\begin{equation}
\label{i_n}
{\rm GL}(n-1,\mathbb C)\ni x\to i_n(x)=x+E_{nn}\in {\rm GL}(n,\mathbb C).
\end{equation}
Then we get
\begin{equation}
\label{Bur-invar}
J_ni_n(\rho_n^{(t)}(\sigma_k))J_n^{-1}=\rho_{n+1}^{(t)}(\sigma_k),\quad 2\leq k\leq n.
\end{equation}

\subsection{The Krammer representation of the braid group $B_n$, different notations}
The {\it Lawrence-Krammer representation} $K_n^{(t,q)}$ of the Braid group $B_n$ was introduced by Lawrence in \cite{Law90}
and studied by Krammer in \cite{Kra00,Kra02} and Bigelow in \cite{Big01,Big02}. To give formulas, we follow Bigelow \cite{Big01} and the corrected version in \cite{Big02}. $K_n^{(t,q)}$ is the following action of $B_n$ on a free
module $V$ of rank
%${n\choose 2}$
$\left(\begin{smallmatrix}n\\2\end{smallmatrix}\right)=\frac{n!}{2!(n-2)!}$
with basis $\{F_{i,j}:1\leq i<j\leq n\}$:
\begin{equation}
\label{Kram-Big}
\sigma_i(F_{j,k})=
\left\{\begin{array}{ll}
F_{j,k}&i\not\in\{j-1,j,k-1,k\},\\
qF_{i,k}+q(q-1)F_{i,j}+(1-q)F_{j,k}&i=j-1,\\
F_{j+1,k}&i=j\not=k-1,\\
qF_{j,i}+(1-q)F_{j,k}+q(1-q)tF_{i,k}&i=k-1\not=j,\\
F_{j,k+1}&i=k,\\
-tq^2F_{j,k}&i=j=k-1.\\
\end{array}\right.
\end{equation}
From \cite{Big02}: ``This action is given in \cite{Kra00}, except that Krammer's $-t$ is my $t$. It is also
in \cite{Big01}, but with a sign error. The name ``Krammer representation'' was
chosen because Krammer seems to have initially found this independently of
Lawrence and without any use of homology''.

We would like to change notations and define the {\it Krammer representation} $k_n^{(t,q)}$ as follows:
\begin{equation}
\label{Kram-Big-new}
\sigma_i(F_{j,k})=
\left\{\begin{array}{ll}
F_{j,k}&i\not\in\{j-1,j,k-1,k\},\\
tF_{i,k}+t(t-1)F_{i,j}+(1-t)F_{j,k}&i=j-1,\\
F_{j+1,k}&i=j\not=k-1,\\
tF_{j,i}+(1-t)F_{j,k}+t(t-1)qF_{i,k}&i=k-1\not=j,\\
F_{j,k+1}&i=k,\\
qt^2F_{j,k}&i=j=k-1.\\
\end{array}\right.
\end{equation}
Connection between two representations is as follows:
\begin{equation}
\label{K-K(new)}
K_n^{(-q,t)}(\sigma_r)=k_n^{(t,q)}(\sigma_r),\quad 1\leq r\leq n-1.
\end{equation}
We think that this choice of the parameters reflects better the situation since, historically Burau used parameter $t$, and $q$, in our notations, is really {\it a parameter of quantization} (see (\ref{K=S^2(B)_q})).
\subsection{Representations of $B_{n+1}$ via representations of the Lie algebra ${\mathfrak gl}_{n}$}
We can generalize the reduced Burau representation of $B_n$  just observing (see \cite{Kos07q}) that ${\rho}_n^{(t)}$ is a product of several  $\exp's$ of an appropriate {\it Serre generators}  of the {\it natural {\rm(}or fundamental {\rm)} representation} of the Lie algebra ${\mathfrak sl}_{n}$ (in fact of ${\mathfrak gl}_{n}$ ) see (\ref{Bur-n+1(pi)}). To be more precise, we show that the reduced Burau representation for $B_3$
\begin{equation}
\label{B_3}
 \sigma_1\mapsto\left(\begin{smallmatrix}
-t&t\\
0&1\\
\end{smallmatrix}\right),\quad
\sigma_2\mapsto \left(\begin{smallmatrix}
1&0\\
1&-t\\
\end{smallmatrix}\right),
\end{equation}
can be rewritten  as follows:
\begin{equation}
\label{B_3-1}
 \sigma_1\mapsto\left(\begin{smallmatrix}
-t&t\\
0&1\\
\end{smallmatrix}\right)=
\left(\begin{smallmatrix}
-t&0\\
0&1\\
\end{smallmatrix}\right)
\left(\begin{smallmatrix}
1&-1\\
0&1\\
\end{smallmatrix}\right)
=\exp(sE_{11})\exp(-E_{12}),
\end{equation}
\begin{equation}
\label{B_3-2}
\sigma_2\mapsto \left(\begin{smallmatrix}
1&0\\
1&-t\\
\end{smallmatrix}\right)=
\left(\begin{smallmatrix}
1&0\\
1&1\\
\end{smallmatrix}\right)\left(\begin{smallmatrix}
1&0\\
0&-t\\
\end{smallmatrix}\right)=
\exp(E_{21})\exp(sE_{22}),
\end{equation}
where $E_{kn}$ are matrix unities and $s=\ln (-t)$ for negative $t$.

For $n=4$ the reduced Burau represenation is defined as follows:
$${\rho}_4^{(t)}:B_4\rightarrow {\rm GL}_3({\mathbb Z}[t,t^{-1}]),$$
\begin{equation}
\label{B_4}
\sigma_1\mapsto\left(\begin{smallmatrix}
-t&t&0\\
0&1&0\\
0&0&1\\
\end{smallmatrix}\right),\,\,
\sigma_2\mapsto\left(\begin{smallmatrix}
1&0&0\\
1&-t&t\\
0&0&1\\
\end{smallmatrix}\right),\,\,
\sigma_3\mapsto\left(\begin{smallmatrix}
1&0&0\\
0&1&0\\
0&1&-t\\
\end{smallmatrix}\right),
\end{equation}
$$
\sigma_1\mapsto
\left(\begin{smallmatrix}
-t&t&0\\
0&1&0\\
0&0&1\\
\end{smallmatrix}\right)=
\left(\begin{smallmatrix}
-t&0&0\\
0&1&0\\
0&0&1\\
\end{smallmatrix}\right)
\left(\begin{smallmatrix}
1&-1&0\\
0&1&0\\
0&0&1\\
\end{smallmatrix}\right)
=\exp(sE_{11})\exp(-E_{12}),
$$
$$
\sigma_2\!\mapsto\!\left(\begin{smallmatrix}
1&0&0\\
1&-t&t\\
0&0&1\\
\end{smallmatrix}\right)\!=\!
\left(\begin{smallmatrix}
1&0&0\\
1&1&0\\
0&0&1\\
\end{smallmatrix}\right)\!
\left(\begin{smallmatrix}
1&0&0\\
0&-t&0\\
0&0&1\\
\end{smallmatrix}\right)
\left(\begin{smallmatrix}
1&0&0\\
0&1&-1\\
0&0&1\\
\end{smallmatrix}\right)
=\exp(E_{21})\exp(sE_{22})\exp(-E_{23}),
$$
$$
\sigma_3\mapsto\left(\begin{smallmatrix}
1&0&0\\
0&1&0\\
0&1&-t\\
\end{smallmatrix}\right)=
\left(\begin{smallmatrix}
1&0&0\\
0&1&0\\
0&1&1\\
\end{smallmatrix}\right)
\left(\begin{smallmatrix}
1&0&0\\
0&1&0\\
0&0&-t\\
\end{smallmatrix}\right)
=\exp(E_{32})\exp(sE_{33}),
$$
where $s=\ln (-t)$ when $t<0$. Finally, we get
\begin{eqnarray}
\label{Bur-exp(sl)3}
\sigma_1=\exp(sE_{11})\exp(-X_1),\,\,
\sigma_2=\exp(Y_1)\exp(sE_{22})\exp(-X_2),\\
\sigma_3=
\exp(Y_2)\exp(sE_{33}).
\end{eqnarray}
In the general case, the formulas are similar, see (\ref{Bur-n+1(pi)}). Below we give the {\it Cartan matrix} $A$
corresponding to a Lie algebra ${\mathfrak sl}_{n}$ and ``the rule to obtain''
the right formulas, where $X_k=E_{kk+1},\,\,Y_k=E_{k+1k}$ and $H_k=E_{kk}-E_{k+1,k+1},$ $1\leq k\leq n-1$, are {\it the Serre generators of ${\mathfrak sl}_{n}$}:
$$
A=\left(\begin{smallmatrix}
2&-1&0&...&0&0\\
-1&2&-1&...&0&0\\
0&-1&2&...&0&0\\
&&&...&&\\
0&0&0&...&2&-1\\
0&0&0&...&-1&2\\
\end{smallmatrix}\right),\quad
\left(\begin{smallmatrix}
E_{11}  &-X_1&&&&&&&&&\\
Y_1&E_{22} &-X_2&&&&&&&&\\
&Y_2&E_{33} &-X_3&&&&&&&\\
&&Y_3&&&&&&&&\\
&&&&&&...&&&&\\
&&&&&&&&&-X_{n-2}&\\
&&&&&&&&&E_{n-1,n-1} &-X_{n-1}\\
&&&&&&&&&Y_{n-1}&E_{nn} \\
\end{smallmatrix}\right).
$$
\begin{lem}
\label{l.Bur-sl_n}
 Let $\pi:{\mathfrak gl}_{n}\to {\rm End}(V)$ be an arbitrary representation
of the
Lie algebra ${\mathfrak gl}_{n}$ in the space $V$. Then the following formulas
\begin{equation}
\begin{array}{l}
\label{Bur-n+1(pi)}
\rho_{n+1}^\pi(\sigma_1)\mapsto\exp(s\pi(E_{11}))\exp(-\pi(X_1)),\\
\rho_{n+1}^\pi(\sigma_k)\mapsto\exp(\pi(Y_{k-1})\exp(s\pi(E_{kk}))\exp(-\pi(X_k)),\\
\rho_{n+1}^\pi(\sigma_n)\mapsto \exp(\pi(Y_{n-1}))\exp(s\pi(E_{nn}))
\end{array}
\end{equation}
gives us  representation $\rho_{n+1}^\pi:B_{n+1}\to {\rm GL}(V)$
of the braid group $B_{n+1}$ in the space $V$. For $\pi$ being {\rm the natural representation} of the Lie algebra ${\mathfrak sl}_{n}$  we
get the reduced Burau representation ${\rho}_{n+1}^{(t)}$.
\end{lem}
%
%\begin{pf} 
It is sufficient to use the fact that formulas (\ref{r-Bur_n}) and (\ref{r-Bur_n(b)}) define representation of the group $B_{n}$.
%\qed
%\end{pf}

\subsection{Connection between the Krammer and the reduced Burau representation}
\begin{thm}
\label{Kr=S^2_q(Bur)}
The Krammer representation $k_n^{(t,q)}$ defined by (\ref{Kram-Big-new}) is equivalent with the {\rm quantization} of the symmetric square of the reduced Burau representation:
\begin{equation}
\label{K=S^2(B)_q}
k_n^{(t,q)}=\big[S^2(\rho_n^{(t)})\big]_q.
\end{equation}
\end{thm}
%%
%\begin{pf}
%
The spectrum of the reduced Burau representation is as follows:\\
${\rm Sp}\,\rho_n^{(t)}(\sigma_{k})=\{-t,\underset{n-1\,\,\text{times}}{1,\dots,1}\}$.
The spectrum of its symmetric square is:
\begin{equation}
\label{Sp(S^2(b_n))}
{\rm Sp}\,S^2(\rho_n^{(t)}(\sigma_{k}))=
\{t^2,
\underset{n-2\,\,\text{times}}{-t\,\,,\dots,\,\,-t,\,\,}
\underset{(n-1)(n-2)/2\,\,\text{times}}{1,\dots,1}\}.
\end{equation}
In \cite[Claim 5.8]{Smel03} (see also \cite[Lemma 5]{Sto09}) the following fact is proved
\begin{lem}
\label{l.Sp(Kr)}
The spectrum of the Krammer representation $k_n^{(t,q)}$ of the group $B_n$, i.e., ${\rm Sp}(\sigma_1^K)$  is as follows:
\begin{equation}
\label{Sp(Kr)}
{\rm Sp}(\sigma_1^K)=\{qt^2,
\underset{n-2\,\,\text{times}}{-t,\dots,-t,}
\underset{(n-1)(n-2)/2\,\,\text{times}}{1,\dots,1}\}.
\end{equation}
\end{lem}
If we set $q=1$  in  (\ref{Sp(Kr)}) we get (\ref{Sp(S^2(b_n))}).
Moreover, in \cite[Proposition 3.2]{Kra00}
the following connection between the Krammer $V$ module and the Burau $W$ module  was proved:
\begin{prop}
\label{pro.Kra}
 If $t = 1$ and $2$ is invertible in $R$, then there is an isomorphism
of $B_n$-modules $V\to S^2W$ {\rm(}symmetric square of $W${\rm)} given by $v_{ij} \mapsto w_{ij}^2$
and, more generally, $v(T ) \mapsto w(T )^2$.
\end{prop}

\begin{rem}
In fact, Krammer showed that the Lawrence-Krammer representation
is a {\it deformation} of the symmetric square of the Burau representation.
We show more, i.e., that this deformation is in fact  {\it quantization}.
\end{rem}

\subsection{Lie algebra $\mathfrak{sl}_2$, its
universal
enveloping algebra $U(\mathfrak{sl}_2)$ and $U(\mathfrak{sl}_2)$-modules}

Recall that the Lie algebra ${\mathfrak sl}_2$ is the  Lie algebra  of $2\times 2$ real matrices with trace equals to zero.
The standard basis $X,\,Y,\,H$ in a Lie algebra ${\mathfrak sl}_2$ is as follows:
\begin{equation}
X= \left(\begin{smallmatrix}
0&1\\
0&0
\end{smallmatrix}\right),\quad Y=
\left(\begin{smallmatrix}
0&0\\
1&0
\end{smallmatrix}\right),\quad H=\left(\begin{smallmatrix}
1&0\\
0&-1
\end{smallmatrix}\right).
\end{equation}
This natural representations $\pi_2:\mathfrak{sl}_2\rightarrow {\rm End}({\mathbb C}^2)$ is called {\it fundamental representation} of the Lie algebra ${\mathfrak sl}_2$.
%%%

All finite-dimensional $U-$module $V$ being the highest weight module of highest weight
$\lambda$ are of the following form (see Kassel,
%1995)
\cite[Theorem V.4.4.]{Kas95})
\begin{equation}
\label{U-mod-1}
\rho(n)(X)\!=\!
S^n(\pi_2(X))\!=\!
\left(\!\begin{smallmatrix}
0&n&0&...&0\\
0&0&n-1&...&0\\
&&&...&\\
0&0&0&...&1\\
0&0&0&...&0
\end{smallmatrix}\!\right),\quad
 \rho(n)(Y)\!=\!
 S^n(\pi_2(Y))\!=\!
 \left(\!\begin{smallmatrix}
0&0  &...&0  &0\\
1&0  &...&0  &0\\
0&2&...&0  &0\\
 &   &...&   &\\
0&0  &...&n&0
\end{smallmatrix}\!\right),\,
\end{equation}
\begin{equation}
\label{U-mod-2}
\rho(n)(H)=
S^n(\pi_2(H))=
\left(\begin{smallmatrix}
n&0      &...&0       &0\\
0  &n-2&...&0       &0\\
   &       &...&        &\\
   &       &...&-n+2&0\\
0  &0      &...&0       &-n\\
\end{smallmatrix}\right),
\end{equation}
where $\lambda ={\rm dim}(V)-1\in{\mathbb N}.$

 \subsection{Pascal triangle  and representations of $B_3$ and Lie algebra $\mathfrak{sl}_2$}
\label{s.2.6}
The connection of the Pascal triangle with representations of $B_3$ was firstly noticed by Humphry
in \cite{Hum00}. Set in (\ref{Rep(q)}) $\Lambda_n=I$ and $q=1$, then
\begin{equation}
\sigma_1\mapsto \sigma_1(1,n),\quad \sigma_2\mapsto \sigma_2(1,n),\quad
\text{where}\quad\sigma_2(1,n)=(\sigma_1(1,n)^{-1})^\sharp
\end{equation}
is an irreducible representations of the group $B_3$ in ${\mathbb C}^{n+1}$. Here the {\it Pascal triangle} denoted by
$\sigma_1(1,n)$, for $n-1\in {\mathbb N}$, is defined as follows:
\begin{equation}
\sigma_1(1,1)=\left(\begin{smallmatrix}
1&1\\
0&1\\
\end{smallmatrix}\right),\quad
\sigma_1(1,2)=\left(\begin{smallmatrix}
1&2&1\\
0&1&1\\
0&0&1\\
\end{smallmatrix}\right),\,\,
\sigma_1(1,3)=\left(\begin{smallmatrix}
1&3&3&1\\
0&1&2&1\\
0&0&1&1\\
0&0&0&1\\
\end{smallmatrix}\right).
\end{equation}
%and so forthe.
\begin{thm}
\label{B_3-sl_2}
Let $\pi_2$
be the fundamental representation of
$\mathfrak{sl}_2$ then
\begin{align}
\sigma_1(1,n)&=\exp(\rho(X))=\exp(S^n(\pi_2(X))),
\\
 \sigma_2(1,n)&=\exp(\rho(-Y))=\exp(S^n(\pi_2(-Y))),
\end{align}
i.e., the Humphry representation of $B_3$ in dimension $n+1$ is the $n^{th}$ symmetric power of the Burau representation
$\rho_3^{-1}$.
\end{thm}
%%%%%%% 
 \subsection{$q$-Pascal triangle  and representations of $B_3$}
\label{s.4.7}
We define  the family of representations of the group $B_3$ using $q$-{\it Pascal triangle} \cite{Kos07q}.
%{KosAlb07q}. 
For $q\in{\mathbb C}^\times\!:={\mathbb C}\setminus\{0\}$
and a matrix
\begin{equation}
\label{l(r)l(n-r)=c}
\Lambda_n={\rm diag}(\lambda_r)_{r=0}^n\quad\text{ such that} \quad
\lambda_r\lambda_{n-r}=const
\end{equation}
 set
\begin{equation}
\label{Rep(q)}
\sigma_1^\Lambda(q,n):=\sigma_1(q,n)D^\sharp_n(q)\Lambda_n,\quad
\sigma_2^\Lambda(q,n):= \Lambda_n^\sharp
D_n(q)\sigma_2(q,n),
\end{equation}
\begin{equation}
\label{D_n(q)} D_n(q)={\rm
diag}(q_r)_{r=0}^n,\,\,q_r=q^{\frac{(r-1)r}{2}}.
\end{equation}
Here  {\it the central symmetry} $A\mapsto A^\sharp$, for the matrix $A\in{\rm Mat}(n+1,{\mathbb C})$, $A=(a_{km})_{0\leq
k,m\leq n}$ is defined as follows:
%$A^\sharp=(a_{km}^\sharp)_{0\leq k,m\leq n},\,\,a_{km}^\sharp=a_{n-k,n-m}$
\begin{equation}
\label{A^sharp=}
\quad\quad
A^\sharp=(a_{km}^\sharp)_{0\leq k,m\leq n},\,\,a_{km}^\sharp=a_{n-k,n-m},
\end{equation}
and $\sigma_1(q,n)$ is the $q$-{\it Pascal triangle} obtained from the Pascal triangle by replacing natural numbers $n$ by  $q$-natural numbers $(n)_q$ (see (\ref{(),[]_q})). More precisely,  define $\sigma_1(q,n)$ as follows $\sigma_1(q,n)=(\sigma_1(q,n)_{km})_{0\leq
k,m\leq n}$, where
\begin{equation}
\label{si_1(q)} \sigma_1(q,n)_{km}=
\left(\begin{smallmatrix} n-k\\
n-m
\end{smallmatrix}\right)_q,\,\,0\leq k,m\leq n,\quad
%\sigma_2(q,n):=(\sigma_1^{-1}(q^{-1},n))^\sharp,
\end{equation}
\begin{equation}
%\label{si_2(q)}
 \sigma_2(q,n):=\sigma_1^{-1}(q^{-1},n)^\sharp,\,\,
\end{equation}
\begin{equation}
\label{si_2(q)} (\sigma_1^{-1}(q^{-1},n))^\sharp_{km}=
\left\{\begin{array}{cc}
0,&\text{\,if\,}\,\,0< k< m\leq n\\
 (-1)^{k+m}q_{k-m}^{-1}C_{k}^{m}(q^{-1}),&\text{\,if\,}\,\,0\leq m\leq k\leq n
\end{array}\right..
\end{equation}
%%%%%%%%%%%%%%%%
Here the $q$-{\it binomial coefficients} or  {\it Gaussian
polynomials} are defined as follows:
\begin{equation}
\label{GP} C_n^k(q):=\left(\begin{smallmatrix} n\\
k
\end{smallmatrix}\right)_q:=\frac{(n)!_q}{(k)!_q(n-k)!_q},\quad
C_n^k[q]:=\left[\begin{smallmatrix} n\\
k
\end{smallmatrix}\right]_q:=\frac{[n]!_q}{[k]!_q[n-k]!_q}
\end{equation}
corresponding to two forms of $q-${\it natural numbers}, defined
by
\begin{equation}
\label{(),[]_q} (n)_q:=\frac{q^n-1}{q-1},\quad
[n]_q:=\frac{q^n-q^{-n}}{q-q^{-1}}.
\end{equation}
%%%%
Define {\it $q$-Pochhammer symbol}
\begin{equation}
\label{(1+x)^k_q}
(a;q)_{n}=\prod _{{k=0}}^{{n-1}}(1-aq^{k})=(1-a)(1-aq)(1-aq^{2})\cdots (1-aq^{{n-1}}).
\end{equation}
We have (see \cite{And76})
\begin{equation}\label{GP2}
(1+x)^k_q:=(-x;q)_{k}=\sum_{r=0}^kq^{r(r-1)/2}C_k^r(q)x^r
=\sum_{r=0}^kq^{r(r-1)/2}
\left(\begin{smallmatrix} k\\
r
\end{smallmatrix}\right)_qx^r.
\end{equation}
%%%%%%%%%%%%%%%%
%Theorem 3~\cite{KosAlb07q,Kos07q}
\begin{thm}[\cite{KosAlb07q,Kos07q}]
\label{t.Rep(q)}
Formulas $  \sigma_1\mapsto
\sigma_1^\Lambda(q,n),\quad \sigma_2\mapsto\sigma_2^\Lambda(q,n)$ {\rm(}see  {\rm(}\ref{Rep(q)}{\rm))} define family of representation of  the group $B_3$ in dimension $n+1$.
\end{thm}
%%%
%%%%%%%%%%%
%
\begin{rem}
\label{pi^Lambda_3}
Let $\pi:B_{3}\to {\rm GL}(m,\mathbb C)$ be a representation. Fix two matrices
$\Lambda_1,\,\Lambda_2\in {\rm GL}(m,\mathbb C)$ and define $\pi^\Lambda(\sigma_1)=\Lambda_1\pi(\sigma_1),\quad
\pi^\Lambda(\sigma_2)=\pi(\sigma_2)\Lambda_2$. When $\pi^\Lambda$ is again a representation of $B_3$?
Under the additional conditions  $\Lambda_1\Lambda_2=cI,\,\,c\in\mathbb C$ and $\Lambda_1={\rm diag}(\lambda_k)_{k=0}^n$ the criteria for $\pi^\Lambda$ to be representation is the following:
\begin{equation}
\Lambda_1\pi(\sigma_1\sigma_2\sigma_1)=\pi(\sigma_1\sigma_2\sigma_1)\Lambda_2.
\end{equation}
This implies condition (\ref{l(r)l(n-r)=c}), i.e.,  $\lambda_r\lambda_{n-r}=c$.
\end{rem}
%%%%%%%%%
\begin{thm} [\cite{Tub01,KosAlb07q,Kos07q}]
%[\cite{Tub01}]
%\cite{KosAlb07q}
\label{Rep(B_3)<5}
All
%irreducible
representations of the group $B_3$ in dimension $\leq 5$ are given by the formulas {\rm(}\ref{Rep(q)}{\rm)}.
\end{thm}
In dimension $2$ we get for $\Lambda_1={\rm diag}(\lambda_0,\lambda_1)$
\begin{equation}
\label{TW(32)(12)}
\sigma_1^{\Lambda_1}(q,1)=
\left(\begin{smallmatrix}
1&1\\
0&1\\
\end{smallmatrix}\right)\left(\begin{smallmatrix}
\lambda_0&0\\
0&\lambda_1\\
\end{smallmatrix}\right),\quad
\sigma_2^{\Lambda_1}(q,1)=
\left(\begin{smallmatrix}
\lambda_1&0\\
0&\lambda_0\\
\end{smallmatrix}\right)\left(\begin{smallmatrix}
1&0\\
-1&1\\
\end{smallmatrix}\right).
\end{equation}
In dimension $3$, we get for $\Lambda_2={\rm diag}(\lambda_0,\lambda_1,\lambda_2)$ with $\lambda_0\lambda_2=\lambda_1^2$
\begin{equation}
\label{TW(33)(12)}
\sigma_1^{\Lambda_2}(q,2)=
\left(\begin{smallmatrix}\lambda_0q&(1+q)\lambda_1&\lambda_2\\
0&\lambda_1&\lambda_2\\
0&0&\lambda_2\\
\end{smallmatrix}\right),\quad
\sigma_2^{\Lambda_2}(q,2)=
 \left(\begin{smallmatrix} \lambda_2&0&0\\
-\lambda_1&\lambda_1&0\\
\lambda_0&-\lambda_0(1+q)&\lambda_0q\\
\end{smallmatrix}\right).
\end{equation}
In dimension $4$, we get for $\Lambda_3={\rm diag}(\lambda_0,\lambda_1,\lambda_2,\lambda_3)$  with $\lambda_0\lambda_3=\lambda_1\lambda_2$:
\begin{equation}
\label{TW(34)1}
 %T_{34}^{\Lambda_4,q}(\sigma_1)
\sigma_1^{\Lambda_3}(q,3)
=\left(\begin{smallmatrix}
1&(1+q+q^2)&(1+q+q^2)&1\\
0&1&(1+q)&1\\
0&0&1&1\\
0&0&0&1
\end{smallmatrix}\right)D_3^\sharp(q)\Lambda_3,
\end{equation}
\begin{equation}
\label{TW(34)2}
\sigma_2^\Lambda(q,3)
=\Lambda_3^\sharp\left(\begin{smallmatrix}
1&0&0&0\\
-1&1&0&0\\
1&-(1+q)&1&0\\
-1&(1+q+q^2)&-(1+q+q^2)&1\\
\end{smallmatrix}\right)D_4(q).
\end{equation}
In dimension $5$, we get for $\Lambda_3={\rm diag}(\lambda_r)_{r=0}^4$
with $\lambda_0\lambda_4=\lambda_1\lambda_3=\lambda_2^2$:
\begin{equation}
\label{TW(35)1}
\sigma_1^\Lambda(q,4)\!=\!
\left(\!\begin{smallmatrix}
1&(1+q)(1+q^2)&(1+q^2)(1+q+q^2)&(1+q)(1+q^2)&1\\
0&1&1+q+q^2&1+q+q^2&1\\
0&0&1&1+q&1\\
0&0&0&1&1\\
0&0&0&0&1
\end{smallmatrix}\!\right)D_4^\sharp(q)\Lambda_4,
\end{equation}
\begin{equation}
\label{TW(35)2}
 \sigma_2^\Lambda(q,4)=
\Lambda^\sharp_4\left(\!\begin{smallmatrix}
1&0&0&0&0\\
-1&1&0&0&0\\
1&-(1+q)&1&0&0\\
-1&(1+q+q^2)&-(1+q+q^2)&1&0\\
1&-(1+q)(1+q^2)&(1+q^2)(1+q+q^2)&-(1+q)(1+q^2)&1\\
\end{smallmatrix}\!\right)D_4(q).
\end{equation}

 \subsection{Quantization of the Pascal triangle}
\begin{rem}
\label{q-Pasc}
The folowing procedure of quantisation (``Quant'') is based on the {\it quantization  of the Pascale triangle} explained in  \cite{KosAlb07q} and 
\cite{Kos07q}, formulas (10) and (11). To be more precise, we change the binomial coefficients
$\left(\begin{smallmatrix} n\\
k
\end{smallmatrix}\right)
$
of the Pascal triangle by  $q$-{\it binomial coefficients}
$
\left(\begin{smallmatrix} n\\
k
\end{smallmatrix}\right)_q\!=\!\frac{(n)!_q}{(k)!_q(n-k)!_q}
$.
\end{rem}
Set $J_n=\sum_{r=0}^nE_{r,n-r}$ and introduce representations  $ T_{3,n+1}^{\Lambda_n,q}$  of the braid group  $B_3$ equivalent to (\ref{Rep(q)}) in dimension $n+1$ as follows:
\begin{equation}
\label{T=TW}
T_{3,n+1}^{\Lambda_n,q,\pi_2}(\sigma_r)=(\sigma_r^\Lambda(q,n))^\sharp=J_n\sigma_r^\Lambda(q,n)J^{-1}_n,\quad r=1,2.
\end{equation}
%%%%%%%%%%%%%%%%%
\begin{rem}
Using \cite{TubWen01}, see also \cite[n7,\,n8]{Kos07q}, we conclude that all representation of $B_3$ in dimension $2,3,4$ and $5$ are equivalent with the following ones, where $\Lambda_1=s\,{\rm diag}(-t,1)$ and $s,t\in {\mathbb C}^\times$:
\begin{equation}
\label{T(32)}
 T_{32}^{\Lambda_1,1,\pi_2}(\sigma_1)\!=\!s\left(\begin{smallmatrix}
-t&t\\
0&1\\
\end{smallmatrix}\right)\!=\!
s\!\left(\begin{smallmatrix}
-t&0\\
0&1\\
\end{smallmatrix}\right)
\left(\begin{smallmatrix}
1&-1\\
0&1\\
\end{smallmatrix}\right)
,\,\,
 T_{32}^{\Lambda_1,1,\pi_2}(\sigma_2)\!=\!s\!\left(\begin{smallmatrix}
1&0\\
1&-t\\
\end{smallmatrix}\right)\!=\!
s\left(\begin{smallmatrix}
1&0\\
1&1\\
\end{smallmatrix}\right)\left(\begin{smallmatrix}
1&0\\
0&-t\\
\end{smallmatrix}\!\right).
\end{equation}
%%%%%%%%%%%%%%%%%
All representation of $B_3$ in dimension $3$  are given by the following  formulas,
where $\Lambda_2={\rm diag}(st^2,-st,s)=s{\rm diag}(t^2,-t,1)$ and $s,t,q\in {\mathbb C}^\times$:
\begin{equation}
\label{T(33)1}
 T_{33}^{\Lambda_2,q,\pi_2}(\sigma_1)=
 %\!\sigma_1^{\Lambda}(q,2)\!=
 \!s\left(\begin{smallmatrix}
t^2q&-t^2(1+q)&t^2\\
0&-t&t\\
0&0&1\\
\end{smallmatrix}\right)=
\!s\left(\begin{smallmatrix}
t^2&0&0\\
0&-t&0\\
0&0&1\\
\end{smallmatrix}\right)
\left(\begin{smallmatrix}
1&-(1+q)&1\\
0&1&-1\\
0&0&1\\
\end{smallmatrix}\right)
\left(\begin{smallmatrix}
q&0&0\\
0&1&0\\
0&0&1\\
\end{smallmatrix}\right),
\end{equation}
\begin{equation}
\label{T(33)2}
 T_{33}^{\Lambda_2,q,\pi_2}(\sigma_2)=
s\left(\begin{smallmatrix}
1&0&0\\
1&-t&0\\
1&-t(1+q)&t^2q\\
\end{smallmatrix}\right)=
s\left(\begin{smallmatrix}
1&0&0\\
1&1&0\\
1&(1+q)&1\\
\end{smallmatrix}\right)
\left(\begin{smallmatrix}
1&0&0\\
0&1&0\\
0&0&q\\
\end{smallmatrix}\right)\left(\begin{smallmatrix}
1&0&0\\
0&-t&0\\
0&0&t^2\\
\end{smallmatrix}\right).
\end{equation}
\begin{rem}
\label{L-K(B_3):)}
If we are sufficiently attentive, we can recognize here the Lawrence-Krammer representation of $B_3$.
\end{rem}

All representation of $B_3$ in dimension $4$ are given by the following  formulas, where
$\Lambda_3={s\rm diag}$ $(-t^3,ut^2,-u^{-1}t,t)$:
\begin{equation}
\label{T(34)1}
 T_{34}^{\Lambda_3,q}(\sigma_1)
=\Lambda_3\left(\begin{smallmatrix}
1&-(1+q+q^2)&(1+q+q^2)&-1\\
0&1&-(1+q)&1\\
0&0&1&-1\\
0&0&0&1
\end{smallmatrix}\right)D_3^\sharp(q),
\end{equation}
\begin{equation}
\label{T(34)2}
 T_{34}^{\Lambda_3,q}(\sigma_2)
%\sigma_1^\Lambda(q,3)
=\left(\begin{smallmatrix}
1&0&0&0\\
1&0&0&0\\
1&(1+q)&1&0\\
1&(1+q+q^2)&(1+q+q^2)&1\\
\end{smallmatrix}\right)D_3(q)\Lambda_3^\sharp.
\end{equation}
For $n=5$ we get
\begin{equation}
\label{T(35)1}
 T_{35}^{\Lambda_4,q}(\sigma_1)=
%\sigma_1^\Lambda(q,4)\!=\!
\Lambda_4\left(\!\begin{smallmatrix}
1&-(1+q)(1+q^2)&(1+q^2)(1+q+q^2)&-(1+q)(1+q^2)&-1\\
0&1&-(1+q+q^2)&(1+q+q^2)&-1\\
0&0&1&-(1+q)&1\\
0&0&0&1&-1\\
0&0&0&0&1
\end{smallmatrix}\!\right)D_4^\sharp(q),
\end{equation}
\begin{equation}
\label{T(35)2}
 T_{35}^{\Lambda_4,q}(\sigma_2)=
 %\sigma_2^\Lambda(q,4)=
%\sigma_1^\Lambda(q,4)\!=\!
\left(\!\begin{smallmatrix}
1&0&0&0&0\\
1&1&0&0&0\\
1&(1+q)&1&0&0\\
1&(1+q+q^2)&(1+q+q^2)&1&0\\
1&(1+q)(1+q^2)&(1+q^2)(1+q+q^2)&(1+q)(1+q^2)&1\\
\end{smallmatrix}\!\right)D_4(q)\Lambda^\sharp_4.
\end{equation}
The general formulas are
\begin{equation}
\label{T(3n)(12)}
 T_{3,n+1}^{\Lambda_{n},q}(\sigma_1)\!=\!\Lambda_{n}D^\sharp_n(q)(\sigma_1^{-1}(q^{-1},n)),\,\,\,
 T_{3,n+1}^{\Lambda_{n},q}(\sigma_2)\!=\!(\sigma_1(q,n))^\sharp D_n(q)\Lambda^\sharp_{n}.
\end{equation}
%%%%%%%%%%%%%%%%%%%%
\end{rem}
\begin{rem}
\label{Irr,n<6}
The irreducibility criteria for the representation (\ref{Rep(q)}) of $B_3$ in dimensions $n\leq 5$ were obtained
in \cite{Tub01,TubWen01}.
\end{rem}
\begin{rem}
\label{dim=6,Irr-?}
For dimensions $n\geq 6$ formulas (\ref{Rep(q)}) gives us some family of  representations of $B_3$ but not all, see
 \cite[Remark 2.11.3]{TubWen01}. For more representations of $B_3$ see, e.g., Westbury, \cite{Westbury95}.
The question of the irreducibility of  representations (\ref{Rep(q)}) is open for $n\geq 6$. 
%Some {\it sufficient conditions}
%of the irreducibility
%are given in \cite{KosAlb07q}.
%
%Representations of $B_3$ defined by (\ref{Rep(q)}) does not exhost all irreducible representations of $B_3$ in dimension $\geq 6$.
%We believe that these conditions are also the {\it necessary} ones.
\end{rem}
%%%%%%%%%%%%%%%%%%%%%%%%%%%%%

\subsection{Symmetric basis in $V\otimes V$}

Let $V$ be a linear space with the basis $(e_k)_{k=1}^n$. The {\it tensor product} $V\otimes V$ is a linear space generated by the basis $(e_k\otimes e_r)_{k,r=1}^n$.
The {\it symmetric square} $S^2(V)$ of the space $V$ is a subspace of $V\otimes V$ generated by one of the equivalent {\it symmetric basis} $e=(e_{kn}^s)_{kn}$ or $v=(v_{kn})_{kn}$,  defined as follows:
\begin{equation}
\label{symm-bas(e)n}
e_{kk}^s=e_k\otimes e_k,\,\,1\leq k\leq n,\quad e_{kr}^s=e_k\otimes e_r+e_r\otimes e_k,\,\,1\leq k<r\leq n,
\end{equation}
\begin{equation}
\label{symm-bas(v)n}
v_{kr}=(e_k+\dots+e_{r-1}) \otimes(e_k+\dots+e_{r-1}),\quad 1\leq k< r\leq n.
\end{equation}
We fix a {\it lexicographic order} on the set $(k,r)$ with $1\leq k\leq r\leq n$. Let an operator $A$ acts on a space $V$. Denote by $A\otimes A$ its tensor product on a space $V\otimes V$ and by $ S^2(A)$ its {\it symmetric square}.
%%%%%%%%%%%%%%%
%%%%%%%%%%%%5
\subsection{Proof of  Theorem~\ref{Kr=S^2_q(Bur)}}
We show that $[S^2(\rho_n^{(t)})_e]_q$, the quantization of the symmetric square $S^2(\rho_n^{(t)})$ of the Burau representations $\rho_n^{(t)}$
for the group $B_n$, calculated in the basis $v$, coincide  with the Lawrence-Krammer representation $k_n^{(t,q)}$ in notations (\ref{Kram-Big-new}), i.e., we show that
\begin{equation}
\label{K-eq-B(n)}
C_n[S^2(\rho^{(t)}_n(\sigma_r))_e]_qC_n^{-1}=k_n^{(t,q)}(\sigma_r),\quad
\text{for}\quad 1\leq r\leq n-1,
\end{equation}
where $C_n$ is the {\it change-of-basis matrix} from the basis $e$  defined by  (\ref{symm-bas(e)n}) to the basis $v$ defined by  (\ref{symm-bas(v)n}) for the vector space $S^2({\mathbb C}^{n-1})$.
%https://en.wikipedia.org/wiki/Change_of_basis

%To prove (\ref{K-eq-B(n)}), taking into consideration (\ref{J(inf)}), (\ref{B(inf)-red-Bur})
%and  (\ref{Bur-invar}),  it is sufficient to consider only three cases, $n=3,4$ and $n=5$ for $\sigma_3$.
The Burau representation $\rho_3^{(t)}$ for $B_3$ on the space $V={\mathbb C}^2$
is as follows:
\begin{equation}
\label{Bur(3)}
 \sigma_1\mapsto \left(\begin{smallmatrix}
-t&t\\
0&1\\
\end{smallmatrix}\right),\quad
\sigma_2\mapsto  \left(\begin{smallmatrix}
1&0\\
1&-t\\
\end{smallmatrix}\right).
\end{equation}
We fix the standard basis $e=(e_k)_{k=1}^n$ in ${\mathbb C}^n$, namely set
\begin{equation}
\label{basis(C)}
e_k=(0,\dots,0,1,0,\dots,0),\quad 1\leq k\leq n.
\end{equation}
In the symmetric square $S^2({\mathbb C}^2)$ the basis $e$ and $v$ are as follows:
\begin{equation}
\label{symm-bas(e)1}
e_{11}^s=e_1\otimes e_1,\quad e_{12}^s=e_1\otimes e_2+e_2\otimes e_1,\quad e_{22}^s=e_2\otimes e_2,
\end{equation}
\begin{equation}
\label{symm-bas(v)1}
v_{12}=e_1\otimes e_1,\quad v_{13}^s=(e_1+e_2)\otimes (e_1+e_2),\quad v_{23}=e_2\otimes e_2.
\end{equation}

The {\it quantization} of the symmetric  square $\big[S^2(\sigma_k)\big]_q$ for $k=1,2$ of the Burau representation $\rho_3^{(t)}$  have the following form in  the basis $e$:
$$
\sigma_1 \mapsto
\left(\begin{smallmatrix}
-t&0\\
0&1\\
\end{smallmatrix}\right)
\left(\begin{smallmatrix}
1&-1\\
0&1\\
\end{smallmatrix}\right)
\stackrel{{\rm Sym^2}}{\mapsto}
\left(\begin{smallmatrix}
t^2&0&0\\
0&-t&0\\
0&0&1
\end{smallmatrix}\right)
\left(\begin{smallmatrix}
1&-2&1\\
0&1&-1\\
0&0&1
\end{smallmatrix}\right)
\stackrel{\text{Quant}}{\mapsto}
$$
\begin{equation}
\label{S^2(s)31}
\left(\begin{smallmatrix}
t^2&0&0\\
0&-t&0\\
0&0&1
\end{smallmatrix}\right)
\left(\begin{smallmatrix}
1&-(1+q)&1\\
0&1&-1\\
0&0&1
\end{smallmatrix}\right)
\left(\begin{smallmatrix}
q&0&0\\
0&1&0\\
0&0&1
\end{smallmatrix}\right)
=
\left(\begin{smallmatrix}
t^2q&-t^2(1+q)&t^2\\
0&-t&t\\
0&0&1
\end{smallmatrix}\right).
\end{equation}
We denote
\begin{equation}
\label{D_(k,2)(q)}
D_{1,2}(q)={\rm diag}(q,1,1)\quad \text{and}\quad D_{2,2}(q)={\rm diag}(1,1,q).
\end{equation}
 Similarly, we get
$$
\sigma_2\mapsto
\left(\begin{smallmatrix}
1&0\\
1&1\\
\end{smallmatrix}\right)
\left(\begin{smallmatrix}
1&0\\
0&-t\\
\end{smallmatrix}\right)
\stackrel{{\rm Sym^2}}{\mapsto}
\left(\begin{smallmatrix}
1&0&0\\
1&1&0\\
1&2&1
\end{smallmatrix}\right)
\left(\begin{smallmatrix}
1&0&0\\
0&-t&0\\
0&0&t^2
\end{smallmatrix}\right)
\stackrel{\text{Quant}}{\mapsto}
$$
\begin{equation}
\label{S^2(s)32}
\left(\begin{smallmatrix}
1&0&0\\
1&1&0\\
1&(1+q)&q
\end{smallmatrix}\right)
\left(\begin{smallmatrix}
1&0&0\\
0&1&0\\
0&0&q
\end{smallmatrix}\right)
\left(\begin{smallmatrix}
1&0&0\\
0&-t&0\\
0&0&t^2
\end{smallmatrix}\right)=
\left(\begin{smallmatrix}
1&0&0\\
1&-t&0\\
1&-t(1+q)&t^2q
\end{smallmatrix}\right).
\end{equation}
%see \cite{KosAlb07q} and \cite{Kos07q},  ,  for quantizing of the Pascal triangle.
Finally, we get
\begin{equation}
\label{T(33)(1,2)}
[S^2(\sigma_1)_e]_q:=
\left(\begin{smallmatrix}
t^2q&-t^2(1+q)&t^2\\
0&-t&t\\
0&0&1
\end{smallmatrix}\right)
,\quad
[S^2(\sigma_2)_e]_q:=\left(\begin{smallmatrix}
1&0&0\\
1&-t&0\\
1&-t(1+q)&t^2q
\end{smallmatrix}\right).
\end{equation}
To enumerate the elements of the basis $e$ and $v$  by the same set of indices, set
\begin{equation}
\label{w-base}
w_{ij}=v_{i,j-1}.
\end{equation}
Since the basis $e$ and $w$ (hence $v$) are connected as follows:
\begin{equation}
\label{e=f(v).3}
%%
\iffalse
\begin{smallmatrix}
e_{11}^s&=&v_{12} &       &\\
e_{12}^s&=&-v_{12}&+v_{13}&-v_{23}\\
e_{22}^s&=&       &       &v_{23}
\end{smallmatrix},\quad
\begin{smallmatrix}
v_{12}&=&e_{11}^s &         &\\
v_{13}&=&e_{11}^s &+e_{12}^s&+e_{22}^s\\
v_{23}&=&           &       &e_{22}^s
\end{smallmatrix}
\fi%%%%%%%%%%%%%
\begin{smallmatrix}
e_{11}^s&=&w_{11} &       &\\
e_{12}^s&=&-w_{11}&+w_{12}&-w_{22}\\
e_{22}^s&=&       &       &w_{22}
\end{smallmatrix},\quad
\begin{smallmatrix}
w_{11}&=&e_{11}^s &         &\\
w_{12}&=&e_{11}^s &+e_{12}^s&+e_{22}^s\\
w_{22}&=&           &       &e_{22}^s
\end{smallmatrix},
\end{equation}
%%%%
the  {\it change-of-basis matrix} $C_3$ in the space $S^2({\mathbb C}^2)$ will be the following:
\begin{equation}
\label{C_3}
C_3=\left(\begin{smallmatrix}
1&-1&0\\
0&1&0\\
0&-1&1
\end{smallmatrix}\right),\quad
C_3^{-1}=\left(\begin{smallmatrix}
1&1&0\\
0&1&0\\
0&1&1
\end{smallmatrix}\right).
\end{equation}
%\commA{repalce $1$ with $r$ and $v$ by $e$} 
We show that
\begin{equation}
\label{K-eq-B(3)}
C_3[S^2(\sigma_r)_e]_qC_3^{-1}=k^{(t,q)}(\sigma_r),\quad
\text{for}\quad r=1,2.
\end{equation}
%{\color{blue}
Indeed, we have
\begin{equation}
\label{S^2(s1)v]q}
[S^2(\sigma_1)_v]_q:=C_3\left(\begin{smallmatrix}
t^2q&-t^2(1+q)&t^2\\
0&-t&t\\
0&0&1
\end{smallmatrix}\right)C_3^{-1}=
\left(\begin{smallmatrix}
t^2q&0&t(t-1)\\
0&0&t\\
0&1&1-t\\
\end{smallmatrix}\right),
\end{equation}
\begin{equation}
\label{S^2(s2)v]q}
[S^2(\sigma_2)_v]_q:=C_3\left(\begin{smallmatrix}
1&0&0\\
1&-t&0\\
1&-t(1+q)&t^2q
\end{smallmatrix}\right)C_3^{-1}=
\left(\begin{smallmatrix}
0&t&0\\
1&1-t&0\\
0&qt(t-1)& t^2q\\
\end{smallmatrix}\right).
\end{equation}
%%%%%%%%%%%%%%%%
This coincides with {\it the Krammer representation $k^{(t,q)}_3$ for  $B_3$} in notations (\ref{Kram-Big-new}):
%
%\commA{\tiny $ABA=BAB$ holds, 27.05.22}
%
\begin{equation}
\label{kra-B3}
%\sigma_1\mapsto\sigma_1^{K^{(q,-t)}_3}:=
k^{(t,q)}_3(\sigma_1)=
\left(\begin{smallmatrix}
t^2q&0&t(t-1)\\
0&0&t\\
0&1&1-t\\
\end{smallmatrix}\right),\quad
k^{(t,q)}_3(\sigma_2)=
\left(\begin{smallmatrix}
0&t&0\\
1&1-t&0\\
0&qt(t-1)& t^2q\\
\end{smallmatrix}\right).
\end{equation}
 The quantization of the symmetric square $S^2(\rho_{n+1}^{(t)})$ of the Burau representations $\rho_{n+1}^{(t)}$
for the group $B_{n+1}$ %{\color{magenta}
will be defined as follows:
%
%\commA{\tiny explain how to do  in the general case, 01.06.22}
%Using Remark~\ref{q-Pasc} define representation of $B_{n+1}$ as follows:
\begin{equation}
\label{s(1)e,m}
%T_{3,m+1}^{\Lambda_m,q,\pi_2}(\sigma_1):=
\sigma_1\to[S^2(\sigma_1)]_q=
%S^2(\Lambda_1)
S^2\big(\exp(sE_{11})\big)[S^2\big(\exp(-E_{12})\big)]_qD_{1,n}(q),
\end{equation}
\begin{equation}
\label{s(n)e,m}
%S_{n+1}^{m,\Lambda}
%T_{3,m+1}^{\Lambda_m,q,\pi_2}(\sigma_n):=[S^m(\sigma_n)_v]_q=
\sigma_n\to
[S^2\big(\exp(E_{n-1,n})\big)]_q
D_{n,n}(q)
S^2\big(\exp(sE_{nn})\big),
%S_2(\Lambda_n),
\end{equation}
\begin{equation}
\label{s(k)e,m}
%S_{n+1}^{m,\Lambda}
%T_{3,m+1}^{\Lambda_m,q,\pi_2}(\sigma_k):=[S^m(\sigma_k)_v]_q=
\sigma_k\to
[S^2\big(\exp(E_{k,k-1}\big)]_qD_{k,n}(q)S^2\big(\exp(sE_{kk})\big)[S^2\big(\exp(-E_{kk+1}\big)]_q,
%Y_{k-1}=E_{k-1,k} X_k=E{kk+1}
\end{equation}
% p.26 02.11.16
%%%%%%%%%%%%%%%%%%%%%%
%\commA{\tiny why, the general case if we replace $S^2$ with  $S^m$ for $m\geq 2$, 01.06.22  }
%
%{\color{magenta}
where $D_{k,n}(q)$ is defined for $1\leq r\leq s\leq n$ as follows:
\begin{equation}
\label{D_(k,n)(q)}
D_{k,n}(q)e^s_{r,s}=qe^s_{r,s}\quad \text{for}\quad k=r=s, \quad\text{and}\quad D_{k,n}(q)e^s_{r,s}=e^s_{r,s}.
\end{equation}
 The Burau representation $\rho_4^{(t)}$  for $B_4$ is as follows:
\begin{equation}
\label{Bur(4)}
 \sigma_1\mapsto
 %\Lambda_{1,3}
% \left(\begin{smallmatrix}
%-t&0&0\\
%0&1&0\\
%0&0&1
%\end{smallmatrix}\right)
 \left(\begin{smallmatrix}
-t&t&0\\
0&1&0\\
0&0&1
\end{smallmatrix}\right),\quad
\sigma_2\mapsto
\left(\begin{smallmatrix}
1&0&0\\
1&-t&t\\
0&0&1
\end{smallmatrix}\right)
%\Lambda_{2,3}
\iffalse
%
 \left(\begin{smallmatrix}
1&0&0\\
0&-t&0\\
0&0&1
\end{smallmatrix}\right)
\left(\begin{smallmatrix}
1&0&0\\
1&1&-1\\
0&0&1
\end{smallmatrix}\right)
%
\fi
%
,\quad
 \sigma_3\mapsto
\left(\begin{smallmatrix}
1&0&0\\
0&1&0\\
0&1&-t
\end{smallmatrix}\right).
%
\iffalse
%
\left(\begin{smallmatrix}
1&0&0\\
0&1&0\\
0&0&-t
\end{smallmatrix}\right),
\Lambda_{3,3},
%
\fi
%
\end{equation}
%where
%$$\Lambda_{1,3}={\rm diag}(-t,1,1),\,\,\Lambda_{2,3}={\rm diag}(1,-t,1),\,\,\Lambda_{3,3}={\rm diag}(1,1,-t),\,\,
%$$
The  symmetric square  of the reduced Burau representation $\rho_4^{(t)}$ for $B_4$ is as follows (see (\ref{B_4}) and (\ref{Bur-exp(sl)3})):
 $$
S^2(\sigma_1)\mapsto
\left(\begin{smallmatrix}
t^2&0&0&0&0&0\\
0&-t&0&0&0&0\\
0&0&1&0&0&0\\
0&0&0&-t&0&0\\
0&0&0&0 &1&0\\
0&0&0&0 &0&1\\
\end{smallmatrix}\right)
\left(\begin{smallmatrix}
1&-2&1&0&0&0\\
0&1&-1&0&0&0\\
0&0&1&0&0&0\\
0&0&0&1&-1&0\\
0&0&0&0&1&0\\
0&0&0&0&0&1\\
\end{smallmatrix}\right)
=
%{\color{blue}
\left(\begin{smallmatrix}
t^2&-2t^2  &t^2&0&0&0\\
0&-t &t&0&0&0\\
0&0&1&0&0&0\\
0&0&0&-t&t&0\\
0&0&0&0&1&0\\
0&0&0&0&0&1\\
\end{smallmatrix}\right),
$$
$$
S^2(\sigma_2)\!\!=\!\!
\left(\begin{smallmatrix}
1&0&0&0&0&0\\
1&1&0&0&0&0\\
1&2&1&0&0&0\\
0&0&0&1&0&0\\
0&0&0&1&1&0\\
0&0&0&0&0&1\\
\end{smallmatrix}\right)\!
\left(\begin{smallmatrix}
1&0&0&0&0&0\\
0&-t&0&0&0&0\\
0&0&t^2&0&0&0\\
0&0&0&1&0&0\\
0&0&0&0&-t&0\\
0&0&0&0&0 &1\\
\end{smallmatrix}\right)\!
\left(\begin{smallmatrix}
1&0&0&0&0&0\\
0&1&0&-1&0&0\\
0&0&1&0&-2&1\\
0&0&0&1 &0&0\\
0&0&0&0&1&-1\\
0&0&0&0&0&1\\
\end{smallmatrix}\right)
\!\!=\!\!
%{\color{red}
%\left(\begin{smallmatrix}
%1&0&0&0&0&0\\
%t&-t&0&-t&0&0\\
%t^2&-2t^2&t^2&2t&-2t&1\\
%0&0&0&1 &0&0\\
%0&0&0&t&-t&1\\
%0&0&0&0 &0 &1\\
%\end{smallmatrix}\right)=}
%{\color{blue}
\left(\begin{smallmatrix}
1&0&0&0&0&0\\
1&-t&0&t&0&0\\
1&-2t&t^2&2t&-2t^2&t^2\\
0&0&0&1 &0&0\\
0&0&0&1&-t&t\\
0&0&0&0 &0 &1\\
\end{smallmatrix}\right),
$$
$$
S^2(\sigma_3)\mapsto
\left(\begin{smallmatrix}
1&0&0&0&0&0\\
0&1&0&0&0&0\\
0&0&1&0&0&0\\
0&1&0&1&0&0\\
0&0&1&0&1&0\\
0&0&1&0&2&1\\
\end{smallmatrix}\right)
\left(\begin{smallmatrix}
1&0&0&0&0&0\\
0&1&0&0&0&0\\
0&0&1&0&0&0\\
0&0&0&-t&0&0\\
0&0&0&0 &-t&0\\
0&0&0&0 &0&t^2\\
\end{smallmatrix}\right)
=
%{\color{blue}
\left(\begin{smallmatrix}
1&0&0&0&0&0\\
0&1&0&0&0&0\\
0&0&1&0&0&0\\
0&1&0&-t&0&0\\
0&0&1&0&-t&0\\
0&0&1&0&-2t&t^2\\
\end{smallmatrix}\right).
$$
%%%%%%%%%%%%%%%%%%%%%%
Using formulas (\ref{s(1)e,m}), (\ref{s(n)e,m}) and (\ref{s(k)e,m}) we get
$$
S^2(\sigma_1)\stackrel{\text{Quant}}{\mapsto}
\left(\begin{smallmatrix}
t^2&0&0&0&0&0\\
0&-t&0&0&0&0\\
0&0&1&0&0&0\\
0&0&0&-t&0&0\\
0&0&0&0     &1&0\\
0&0&0&0     &0&1\\
\end{smallmatrix}\right)
\left(\begin{smallmatrix}
q&-(1+q)&1&0&0&0\\
0&1&-1&0&0&0\\
0&0&1&0&0&0\\
0&0&0&1&1&0\\
0&0&0&0&1&0\\
0&0&0&0&0&1\\
\end{smallmatrix}\right)
=
%{\color{blue}
\left(\begin{smallmatrix}
t^2q&-t^2(1+q)&t^2&0&0&0\\
0&-t&t&0&0&0\\
0&0&1&0&0&0\\
0&0&0&-t&t&0\\
0&0&0&0&1&0\\
0&0&0&0&0&1\\
\end{smallmatrix}\right),
$$
$$
S^2(\sigma_3)\!\!\stackrel{\text{Quant}}{\mapsto}\!\!
\left(\begin{smallmatrix}
1&0&0&0&0&0\\
0&1&0&0&0&0\\
0&0&1&0&0&0\\
0&1&0&1&0&0\\
0&0&1&0&1&0\\
0&0&1&0&(1+q)&q\\
\end{smallmatrix}\right)
\left(\begin{smallmatrix}
1&0&0&0&0&0\\
0&1&0&0&0&0\\
0&0&1&0&0&0\\
0&0&0&-t&0&0\\
0&0&0&0     &-t&0\\
0&0&0&0     &0&t^2\\
\end{smallmatrix}\right)\!
=
%{\color{blue}
\left(\begin{smallmatrix}
1&0&0&0&0&0\\
0&1&0&0&0&0\\
0&0&1&0&0&0\\
0&1&0&-t&0&0\\
0&0&1&0&-t&0\\
0&0&1&0&-t(1+q)&t^2q\\
\end{smallmatrix}\right).
$$
$$
S^2(\sigma_2)\!\!\stackrel{\text{Quant}}{\mapsto}\!\!
\left(\begin{smallmatrix}
1&0&0&0&0&0\\
1&1&0&0&0&0\\
1&1+q&1&0&0&0\\
0&0&0&1&0&0\\
0&0&0&1&1&0\\
0&0&0&0&0&1\\
\end{smallmatrix}\right)\!
\left(\begin{smallmatrix}
1&0&0&0&0&0\\
0&-t&0&0&0&0\\
0&0&qt^2&0&0&0\\
0&0&0&1&0&0\\
0&0&0&0&-t&0\\
0&0&0&0&0 &1\\
\end{smallmatrix}\right)\!
\left(\begin{smallmatrix}
1&0&0&0&0&0\\
0&1&0&-1&0&0\\
0&0&1&0&-(1+q)&1\\
0&0&0&1 &0&0\\
0&0&0&0&1&-1\\
0&0&0&0 &0 &1\\
\end{smallmatrix}\right)
\!\!=\!\!
$$
$$
\left(\begin{smallmatrix}
1&0&0&0&0&0\\
1&-t&0&0&0&0\\
1&-t(1+q)&qt^2&0&0&0\\
0&0&0&1&0&0\\
0&0&0&1&-t&0\\
0&0&0&0&0&1\\
\end{smallmatrix}\!\right)\!
\left(\!\begin{smallmatrix}
1&0&0&0&0&0\\
0&1&0&-1&0&0\\
0&0&1&0&-(1+q)&1\\
0&0&0&1 &0&0\\
0&0&0&0&1&-1\\
0&0&0&0 &0 &1\\
\end{smallmatrix}\right)\!=\!
%
%{\color{blue}
\left(\!\begin{smallmatrix}
1&0&0&0&0&0\\
1&-t&0&t&0&0\\
1&-t(1+q)&t^2q&t(1+q)&-t^2(1+q)&t^2\\
0&0&0&1 &0&0\\
0&0&0&1&-t&t\\
0&0&0&0 &0 &1\\
\end{smallmatrix}\!\right).
$$
% p.27 02.11.16
%%%
Finally, we get in the basis $e$:
\begin{equation}\label{s(1,3)e}
%\sigma_1\!\!\stackrel{S^2_Q(\text{Bur})}{\mapsto}\!
%
[S^2(\sigma_1)_e]_q:=
\left(\begin{smallmatrix}
t^2q&-t^2(1+q)&t^2&0&0&0\\
0&-t&t&0&0&0\\
0&0&1&0&0&0\\
0&0&0&-t&t&0\\
0&0&0&0&1&0\\
0&0&0&0&0&1\\
\end{smallmatrix}\right),\quad
%\sigma_3\!\!\stackrel{S^2_Q(\text{Bur})}{\mapsto}\!
[S^2(\sigma_3)_e]_q:=
\left(\begin{smallmatrix}
1&0&0&0&0&0\\
0&1&0&0&0&0\\
0&0&1&0&0&0\\
0&1&0&-t&0&0\\
0&0&1&0&-t&0\\
0&0&1&0&-t(1+q)&t^2q\\
\end{smallmatrix}\right),
\end{equation}
\begin{equation}\label{s(2)e}
%\sigma_2\stackrel{S^2_Q(\text{Bur})}{\mapsto}
[S^2(\sigma_2)_e]_q:=
\left(\begin{smallmatrix}
1&0&0&0&0&0\\
1&-t&0&t&0&0\\
1&-t(1+q)&t^2q&t(1+q)&-t^2(1+q)&t^2\\
0&0&0&1 &0&0\\
0&0&0&1&-t&t\\
0&0&0&0 &0 &1\\
\end{smallmatrix}\right).
\end{equation}
We rewrite (\ref{s(1,3)e}) and  (\ref{s(2)e}) in the basis $v$ and compare with the expressions $k_4^{(t,q)}(\sigma_r),$ for $r\!=\!1,2,3$. As before, we keep the lexicographic order to enumerate the elements of the basis $e$ and $w$ in $S^2({\mathbb C}^3)$:
$$
e=(e_{11}^s,e_{12}^s,e_{22}^s,e_{13}^s,e_{23}^s,e_{33}^s),\quad w=(w_{11},w_{12},w_{22},w_{13},w_{23},w_{33}).
$$
%%%%%%%%%%%%%%%%
\iffalse
The expressions of $e$ in terms of $w$ and vice versa are the following:
\begin{equation}
\label{e=f(v).4}
\begin{smallmatrix}
e_{11}^s&=&w_{11} &       &       &&&\\
e_{12}^s&=&-w_{11}&+w_{12}&-w_{22}&&&\\
e_{22}^s&=&       &       & w_{22}&&&\\
e_{13}^s&=&       &-w_{12}&+ w_{22}&+w_{13}&- w_{23}&\\
e_{23}^s&=&       &       &-w_{22} &       &+ w_{23}&-w_{33}\\
e_{33}^s&=&       &       &        &       &        &w_{33}\\
\end{smallmatrix},
\end{equation}
%
\begin{equation}
\label{v=f(e).4}
\begin{smallmatrix}
w_{11}&=&e_{11}^s &         &&&&\\
w_{12}&=&e_{11}^s &+e_{12}^s&+e_{22}^s&&&\\
w_{22}&=&           &       &e_{22}^s&&&\\
w_{13}&=&e^s_{11}&+e^s_{12}&+ e^s_{22}&+e^s_{13}&+ e^s_{23}&+e^s_{33}\\
w_{23}&=&&&+ e^s_{22}&&+ e^s_{23}&+e^s_{33}\\
w_{33}&=&&&&&&e^s_{33}\\
\end{smallmatrix}.
\end{equation}
\fi
%%%%%%%%%%%%%

%Similarly to (\ref{C_3}),
By Lemma~\ref{l.C(-1)n,C_n},
the   change-of-basis matrices $C_4$ and $C_4^{-1}$ in the space $S^2({\mathbb C}^4)$ are the following:
\begin{equation}
\label{C_4}
C_4=\left(\begin{smallmatrix}
1&-1&0& 0&0&0\\
0& 1&0&-1&0&0\\
0&-1&1& 1&-1&0\\
0&0 &0& 1&0&0\\
0&0 &0&-1& 1&0\\
0&0 &0& 0&-1&1
\end{smallmatrix}\right),\quad
C_4^{-1}=\left(\begin{smallmatrix}
1&1&0&1 &0&0\\
0&1&0&1 &0&0\\
0&1&1&1 &1&0\\
0&0&0&1 &0&0\\
0&0&0&1 &1&0\\
0&0&0&1 &1&1
\end{smallmatrix}\right).
\end{equation}
We calculate $[S^2(\sigma_r)_v]_q=C_4[S^2(\sigma_r)_e]_qC_4^{-1}$ for $r=1,2,3$. We have
%
%\commA{\eqref{s1},
%\eqref{s2},\eqref{s3}}
%
\begin{eqnarray}
\label{s1}
&&
[S^2(\sigma_1)_v]_q=C_4
\left(\begin{smallmatrix}
t^2q&-{\color{blue}t^2}(1+q)&{\color{blue}t^2}&0&0&0\\
0&-t&{\color{blue}t}&0&0&0\\
0&0&1&0&0&0\\
0&0&0&-t&{\color{blue}t}&0\\
0&0&0&0&1&0\\
0&0&0&0&0&1\\
\end{smallmatrix}\right)
C_4^{-1}=
\left(\begin{smallmatrix}
qt^2&0&t(t-1)&0&t(t-1)&0\\
0    &0&t&0&0&0\\
0    &1&1-t&0&0&0\\
0&0&0&0&t&0\\
0&0&0&1&1-t&0\\
0&0&0&0&0&1\\
\end{smallmatrix}\right),\\
%%
%\label{s2}
\nonumber
&&
[S^2(\sigma_2)_v]_q\!=\!C_4\!
\left(\begin{smallmatrix}
1&0&0&0&0&0\\
{\color{blue}1}&-t&0&{\color{blue}t}&0&0\\
1&-t(1+q)&t^2q&t(1+q)&-t^2(1+q)&t^2\\
0&0&0&1 &0&0\\
0&0&0&{\color{blue}1}&-t&{\color{blue}t}\\
0&0&0&0 &0 &1\\
\end{smallmatrix}\right)\!
C_4^{-1}\!=\!
\left(\!\begin{smallmatrix}
0&t&0&0&0&0\\
1&1-t&0&0&0&0\\
0&qt(t-1)&qt^2&0&0&t(t-1)\\
0&0&0&1&0&0\\
0&0&0&0&0&t\\
0&0&0&0&1&1-t\\
\end{smallmatrix}\!\right),\\
%%%%
\label{s3}
&&
[S^2(\sigma_3)_v]_q=C_4
\left(\begin{smallmatrix}
1&0&0&0&0&0\\
0&1&0&0&0&0\\
0&0&1&0&0&0\\
0&1&0&-t&0&0\\
0&0&1&0&-t&0\\
0&0&1&0&-t(1+q)&t^2q\\
\end{smallmatrix}\right)
C_4^{-1}=
\left(\begin{smallmatrix}
1&0&0&0&0&0\\
0&0&0&t&0&0\\
0&0&0&0&t&0\\
0&1&0&1-t&0  &0\\
0&0&1&0  &1-t&0\\
0&0&0&qt(t-1)&qt(t-1)&qt^2\\
\end{smallmatrix}\right).
\end{eqnarray}
% p.27 02.11.16
%
This coincide with {\it the Krammer representation $k_4^{(t,q)}(\sigma_r)$ for  $B_4$} in notations (\ref{Kram-Big-new}):
 \begin{equation*}
%\label{LowKra-B-4}
%\sigma_1\mapsto
k^{(t,q)}_4(\sigma_1)=
\left(\begin{smallmatrix}
qt^2&0&t(t-1)&0&t(t-1)&0\\
0    &0&t&0&0&0\\
0    &1&1-t&0&0&0\\
0&0&0&0&t&0\\
0&0&0&1&1-t&0\\
0&0&0&0&0&1\\
\end{smallmatrix}\right),\quad
%\sigma_2\mapsto
k^{(t,q)}_4(\sigma_2)=
\left(\begin{smallmatrix}
0&t&0&0&0&0\\
1&1-t&0&0&0&0\\
0&qt(1-t)&qt^2&0&0&t(t-1)\\
0&0&0&1&0&0\\
0&0&0&0&0&t\\
0&0&0&0&1&1-t\\
\end{smallmatrix}\right),
\end{equation*}
\begin{equation*}
%\sigma_3\mapsto
k^{(t,q)}_4(\sigma_3)=
\left(\begin{smallmatrix}
1&0&0&0&0&0\\
0&0&0&t&0&0\\
0&0&0&0&t&0\\
0&1&0&1-t&0  &0\\
0&0&1&0  &1-t&0\\
0&0&0&qt(t-1)&
qt(t-1)&qt^2\\
\end{smallmatrix}\right).
\end{equation*}
%%%%%%%%%%
\iffalse
For the general $B_{n+1}$ the formulas $[S^2(\sigma_k)e]_q$ are similar to (\ref{s(1,3)e}) and (\ref{s(2)e}),
%and have the following form,
taking into consideration Lemma~\ref{l.Bur-sl_n} and
%Remark~\ref{q(12)-place}:
\begin{equation}
\label{s(1)e}
[S^2(\sigma_1)_e]_q=S^2\big(\exp(sE_{11})\big)[S^2\big(\exp(-X_1)\big)]_q,
\end{equation}
\begin{equation}
\label{s(n)e}
[S^2(\sigma_n)_e]_q=[S^2\big(\exp(Y_{n-1})\big)]_qS^2\big(\exp(sE_{nn})\big),
\end{equation}
\begin{equation}
\label{s(k)e}
[S^2(\sigma_k)_e]_q=
[S^2\big(\exp(Y_{k-1}\big)]_qS^2\big(\exp(sE_{kk})\big)[S^2\big(\exp(-X_k\big)]_q.
\end{equation}
\fi
%%%%%%%%%%
%{\color{blue}
\begin{rem}
\label{r.Kash}
Rinat Kashaev pointed out to the author 
during his visit to Geneva 
in May 2022, that the formulas \eqref{s1}--\eqref{s3} are not correct.
\end{rem}
%%%%%%%%%%%%%%%%%
%%%%%%%%%%%%%%%%
\begin{lem}
\label{l.C(-1)n,C_n}
The matrices $C_{n+1}^{-1}$ and $C_{n+1}$ are defined as follows.
The columns of the matrix $C_{n+1}^{-1}$ are defined by the following relations:
\begin{equation}
\label{C(-1)n}
\quad\quad\quad\quad
w_{ij}=\sum_{i\leq k\leq r\leq j}e^s_{kr},\quad 1\leq i\leq j\leq n.
\end{equation}
The columns of the matrix $C_{n+1}$ are defined by the following relations:
%compare with (\ref{e=f(v).3}) and (\ref{e=f(v).4}):
\begin{equation}
\label{C_n}
e^s_{ij}=
\left\{\begin{array}{lll}
(-1)^{i+j}\sum_{i\leq k\leq r\leq j}(-1)^{k+r}w_{kr},&\text{if}&0\leq j-i\leq 1,\\
\sum_{r=i+1}^{j-1}
%\sum_{k=i}^{i+1}(-1)^kw_{kr}-\sum_{k=i}^{i+1}(-1)^kw_{kj},
(-w_{ir}+w_{i+1,r})+w_{ij}-w_{i+1,j},&\text{if}&i=1,\,\,\,\,2\leq j-i,\\
-w_{ii}+\sum_{r=i+1}^{j-1}
%\sum_{k=i}^{i+1}(-1)^kw_{kr}-\sum_{k=i}^{i+1}(-1)^kw_{kj},
(-w_{ir}+w_{i+1,r})+w_{ij}-w_{i+1,j},&\text{if}&i>1,\,\,\,\,2\leq j-i.
\end{array}\right.
\end{equation}
\end{lem}
%%%
%\begin{pf}
%
Indeed, to prove (\ref{C_3}), (\ref{C_4}) and (\ref{C_5}) we note that matrix $E_n=C_n^{-1}$ is block upper-triangular
\begin{equation}
\label{E_n^{-1}}
E_n=\left(\begin{smallmatrix}
a_n&b_n\\
0&c_n\\
\end{smallmatrix}\right),\quad\text{therefore}\quad
C_n=E_n^{-1}=\left(\begin{smallmatrix}
a_n^{-1}&-a_n^{-1}b_nc_n^{-1}\\
0&c_n^{-1}\\
\end{smallmatrix}\right).
\end{equation}
For $n=3$ we have
$$
E_3=\left(\begin{smallmatrix}
1&1&0\\
0&1&0\\
0&1&1
\end{smallmatrix}\right)=\left(\begin{smallmatrix}
a_3&b_3\\
0&c_3\\
\end{smallmatrix}\right),\quad
E_3^{-1}=\left(\begin{smallmatrix}
a_3^{-1}&-a_3^{-1}b_3c_3^{-1}\\
0&c_3^{-1}\\
\end{smallmatrix}\right)
$$
where
$$
a_3=(1),\quad b_3=(1\,\,0),\quad c_3=
\left(\begin{smallmatrix}
1&0\\
1&1\\
\end{smallmatrix}\right).
$$
We have $c_3^{-1}=
\left(\begin{smallmatrix}
1&0\\
-1&1\\
\end{smallmatrix}\right)$, therefore $a_3^{-1}b_3c_3^{-1}=(1)(1\,\,0)\left(\begin{smallmatrix}
1&0\\
1&1\\
\end{smallmatrix}\right)=(1\,\,0)$, hence
$
C_3=\left(\begin{smallmatrix}
1&-1&0\\
0&1&0\\
0&-1&1
\end{smallmatrix}\right).
$
For $n=4$ we get
$$
E_4=\left(\begin{smallmatrix}
1&1&0&1 &0&0\\
0&1&0&1 &0&0\\
0&1&1&1 &1&0\\
0&0&0&1 &0&0\\
0&0&0&1 &1&0\\
0&0&0&1 &1&1
\end{smallmatrix}\right)=\left(\begin{smallmatrix}
a_4&b_4\\
0&c_4\\
\end{smallmatrix}\right),\quad
E_4^{-1}=\left(\begin{smallmatrix}
a_4^{-1}&-a_4^{-1}b_4c_4^{-1}\\
0&c_4^{-1}\\
\end{smallmatrix}\right)
$$
where
$$
a_4=E_3,\quad b_4=\left(\begin{smallmatrix}
1&0&0\\
1&0&0\\
1&1&0\\
\end{smallmatrix}\right),
\quad c_4=
\left(\begin{smallmatrix}
1&0&0\\
1&1&0\\
1&1&1\\
\end{smallmatrix}\right).
$$
Since $c_4^{-1}=
\left(\begin{smallmatrix}
1&0&0\\
-1&1&0\\
0&-1&1\\
\end{smallmatrix}\right)$ we get
$$
a_4^{-1}b_4c_4^{-1}=\left(\begin{smallmatrix}
1&-1&0\\
0&1&0\\
0&-1&1
\end{smallmatrix}\right)
\left(\begin{smallmatrix}
1&0&0\\
1&0&0\\
1&1&1
\end{smallmatrix}\right)
\left(\begin{smallmatrix}
1&0&0\\
-1&1&0\\
0&-1&1\\
\end{smallmatrix}\right)=\left(\begin{smallmatrix}
0&0&0\\
1&0&0\\
-1&1&0\\
\end{smallmatrix}\right),\quad
C_4=\left(\begin{smallmatrix}
1&-1&0& 0&0&0\\
0& 1&0&-1&0&0\\
0&-1&1& 1&-1&0\\
0&0 &0& 1&0&0\\
0&0 &0&-1& 1&0\\
0&0 &0& 0&-1&1
\end{smallmatrix}\right).
$$
For $n=5$ we have
\iffalse
\left(\begin{smallmatrix}
1|&x_{12}&x_{13}&x_{14}&|x_{15}\\
\\\hline\\
0|&1&x_{23}&x_{24}&|x_{25}\\
0|&0&1&x_{34}&|x_{35}\\
0|&0&0&1&|x_{45}\\ \hline
0|&0&0&0&|\,\,\,1\,\,\,\\
\end{smallmatrix}\right)
\fi
$$
E_5=C_5^{-1}=\left(\begin{smallmatrix}
1|&1&0|&1 &0&0|&1&0&0&0\\  \hline
0|&1&0|&1 &0&0|&1&0&0&0\\
0|&1&1|&1 &1&0|&1&1&0&0\\ \hline
0|&0&0|&1 &0&0|&1&0&0&0\\
0|&0&0|&1 &1&0|&1&1&0&0\\
0|&0&0|&1 &1&1|&1&1&1&0\\ \hline
0|&0&0|&0&0&0|&1&0&0&0\\
0|&0&0|&0&0&0|&1&1&0&0\\
0|&0&0|&0&0&0|&1&1&1&0\\
0|&0&0|&0&0&0|&1&1&1&1\\
\end{smallmatrix}\right)=\left(\begin{smallmatrix}
a_5&b_5\\
0&c_5\\
\end{smallmatrix}\right),\quad
E_5^{-1}=\left(\begin{smallmatrix}
a_5^{-1}&-a_5^{-1}b_5c_5^{-1}\\
0&c_5^{-1}\\
\end{smallmatrix}\right),
$$
where
$$
a_5=E_4,\quad b_5=\left(\begin{smallmatrix}
1&0&0&0\\
1&0&0&0\\
1&1&0&0\\
1&0&0&0\\
1&1&0&0\\
1&1&1&0
\end{smallmatrix}\right),
\quad c_5=
\left(\begin{smallmatrix}
1&0&0&0\\
1&1&0&0\\
1&1&1&0\\
1&1&1&1\\
\end{smallmatrix}\right).
$$
Since $c_5^{-1}=
\left(\begin{smallmatrix}
1&0&0&0\\
-1&1&0&0\\
0&-1&1&0\\
0&0&-1&1
\end{smallmatrix}\right)$
we get
$$
a_5^{-1}b_5c_5^{-1}=\left(\begin{smallmatrix}
1&-1&0& 0&0&0\\
0& 1&0&-1&0&0\\
0&-1&1& 1&-1&0\\
0&0 &0& 1&0&0\\
0&0 &0&-1& 1&0\\
0&0 &0& 0&-1&1
\end{smallmatrix}\right)
\left(\begin{smallmatrix}
1&0&0&0\\
1&0&0&0\\
1&1&0&0\\
1&0&0&0\\
1&1&0&0\\
1&1&1&0
\end{smallmatrix}\right)
\left(\begin{smallmatrix}
1&0&0&0\\
-1&1&0&0\\
0&-1&1&0\\
0&0&-1&1
\end{smallmatrix}\right)=\left(\begin{smallmatrix}
0&0&0&0\\
1&0&0&0\\
-1&1&0&0\\
1&0&0&0\\
-1&1&0&0\\
0&-1&1&0\\
\end{smallmatrix}\right),
$$
and
$$
\hskip 2 cm
C_5=\left(\begin{smallmatrix}
1|&-1&0|& 0&0&0|   &0 &0 &0&0\\ \hline
0|& 1&0|&-1&0&0|   &-1&0 &0&0\\
0|&-1&1|& 1&-1&0|  & 1&-1&0&0\\\hline
0|&0 &0|& 1&0&0 |  &-1&0 &0&0\\
0|&0 &0|&-1& 1&0|&1&-1&0 &0\\
0|&0 &0|& 0&-1&1|  &0 &1 &-1&0\\\hline
0|&0&0|&0&0&0|& 1& 0& 0&0\\
0|&0&0|&0&0&0|&-1& 1& 0&0\\
0|&0&0|&0&0&0|&0 &-1& 1&0\\
0|&0&0|&0&0&0|&0 & 0&-1&1\\
\end{smallmatrix}\right).\hskip 3.5 cm\Box
$$
%\end{pf}
%
\begin{rem}
\label{r.C_n=}
The matrix $E_n=C_n^{-1}$ and  $E_n^{-1}=C_n$ have the following form
\begin{equation}
\label{C_n=}
E_n=\left(\begin{smallmatrix}
e_{11}&e_{12}&e_{13}&\dots&e_{1n-1}\\
0     &e_{22}&e_{23}&\dots&e_{2n-1}\\
0     &0     &e_{33}&\dots&e_{3n-1}\\
      &      &      &\dots&\\
0     &0     &0     &\dots&e_{n-1,n-1}\\
\end{smallmatrix}\right),\quad
E_n^{-1}=\left(\begin{smallmatrix}
e_{11}^{-1}&e_{12}^{-1}&e_{13}^{-1}&\dots&e_{1n-1}^{-1}\\
0     &e_{22}^{-1}&e_{23}^{-1}&\dots&e_{2n-1}^{-1}\\
0     &0     &e_{33}^{-1}&\dots&e_{3n-1}^{-1}\\
      &      &      &\dots&\\
0     &0     &0     &\dots&e_{n-1,n-1}^{-1}\\
\end{smallmatrix}\right),
\end{equation}
where $e_{kr}\in{\rm Mat}(k\times r, {\mathbb C})$ for $1\leq k\leq r\leq n-1$ are as follows:
$$
e_{11}=(1),\,\,e_{22}=\left(\begin{smallmatrix}
1&0\\
1&1
\end{smallmatrix}\right),\,\,e_{33}=\left(\begin{smallmatrix}
1&0&0\\
1&1&0\\
1&1&1\\
\end{smallmatrix}\right),\,\,\,e_{kk}=(I-e_k)^{-1},\,\,\,e_k=\sum_{r=1}^{k-1}E_{r+1,r},
$$
$e_{kr}=(e_{kk}, 0_{k,r-k})$, where $0_{k,r}=0$ in ${\rm Mat}(k\times r, {\mathbb C})$. For $E_n^{-1}$ we have
$$
e_{11}^{-1}=(1),\,\,e_{22}^{-1}=\left(\begin{smallmatrix}
1&0\\
-1&1
\end{smallmatrix}\right),\,\,e_{33}^{-1}=\left(\begin{smallmatrix}
1&0&0\\
-&1&0\\
0&-1&1\\
\end{smallmatrix}\right),\,\,\,e_{kk}=(I-e_k),
$$
and $e_{1r}^{-1}=0$ for $2\leq r\leq n-1$, $e_{kr}^{-1}:=(-e_{kk}^{-1},0_{k,r-k})$ for $2\leq k<r\leq n-1$.
Notation $e_{kr}^{-1}$ for $k<r$ does not mean the inverse matrix to $e_{kr}$!
\end{rem}
To prove (\ref{K-eq-B(n)}) for general $B_{n
+1}$, it is sufficient to consider only three cases $n=3,4$ and $n=5$ for $\sigma_3$.
Indeed, we can use the {\it invariance} of the reduced Burau representation of $B_\infty$ and $B_{n}$,
(see  (\ref{J(inf)}), (\ref{B(inf)-red-Bur}) and  (\ref{Bur-invar})) and invariance of the Lawrence-Krammer representation (\ref{Kram-Big-new}).

%%%%%
The Burau representation $\rho_5^{(t)}$  for $B_5$ is as follows:
\begin{equation}
\label{Bur(5)}
 \sigma_1\mapsto \left(\begin{smallmatrix}
-t&t&0&0\\
0&1&0&0\\
0&0&1&0\\
0&0&0&1
\end{smallmatrix}\right),\,\,\,
\sigma_2\mapsto  \left(\begin{smallmatrix}
1&0&0&0\\
1&-t&t&0\\
0&0&1&0\\
0&0&0&1
\end{smallmatrix}\right),\,\,\,
 \sigma_3\mapsto
\left(\begin{smallmatrix}
1&0&0&0\\
0&1&0&0\\
0&1&-t&t\\
0&0&0&1
\end{smallmatrix}\right),\,\,\,
 \sigma_4\mapsto
\left(\begin{smallmatrix}
1&0&0&0\\
0&1&0&0\\
0&0&1&0\\
0&0&1&-t
\end{smallmatrix}\right).
\end{equation}
By Lemma~\ref{l.C(-1)n,C_n}, the   change-of-basis matrix $C_5$  and $C_5^{-1}$ in the space $S^2({\mathbb C}^5)$ are the following:
\begin{equation}
\label{C_5}
C_5=\left(\begin{smallmatrix}
1&-1&0& 0&0&0   &0 &0 &0&0\\
0& 1&0&-1&0&0   &-1&0 &0&0\\
0&-1&1& 1&-1&0  & 1&-1&0&0\\
0&0 &0& 1&0&0   &-1&0 &0&0\\
0&0 &0&-1& 1&0&1&-1&0 &0\\
0&0 &0& 0&-1&1  &0 &1 &-1&0\\
0&0&0&0&0&0& 1& 0& 0&0\\
0&0&0&0&0&0&-1& 1& 0&0\\
0&0&0&0&0&0&0 &-1& 1&0\\
0&0&0&0&0&0&0 & 0&-1&1\\
\end{smallmatrix}\right),\quad
C_5^{-1}=\left(\begin{smallmatrix}
1&1&0&1 &0&0&1&0&0&0\\
0&1&0&1 &0&0&1&0&0&0\\
0&1&1&1 &1&0&1&1&0&0\\
0&0&0&1 &0&0&1&0&0&0\\
0&0&0&1 &1&0&1&1&0&0\\
0&0&0&1 &1&1&1&1&1&0\\
0&0&0&0&0&0&1&0&0&0\\
0&0&0&0&0&0&1&1&0&0\\
0&0&0&0&0&0&1&1&1&0\\
0&0&0&0&0&0&1&1&1&1\\
\end{smallmatrix}\right).
\end{equation}
Using  (\ref{Kram-Big-new}) we get
\begin{equation*}
%\sigma_3\mapsto
k^{(t,q)}_5(\sigma_3)=
\left(\begin{smallmatrix}
1&0&0   &0&0&0&0&0&0&0\\
0&0&t   &0&0&0&0&0&0&0\\
0&1&1-t&0&0&0&0&0&0&0\\
0&0&0&1&0&0&0&0&0&0\\
0&0&0&0&0&1&0&0&0&0\\
0&0&0&0&1&1-t&0&0&0&0\\
0&0&0&0&0&0&1&0&0&0\\
0&0&qt(t-1)&0&0&qt(t-1)&0&qt^2&0&t(t-1)\\
0&0&0&0&0&0&0&0&1&t\\
0&0&0&0&0&0&0&0&0&1-t\\
\end{smallmatrix}\right).
\end{equation*}
Indeed, since
$$
\sigma_3\!\mapsto\!\left(\begin{smallmatrix}
1&0&0&0\\
0&1&0&0\\
0&1&-t&t\\
0&0&0&1\\
\end{smallmatrix}\right)\!=\!
\left(\begin{smallmatrix}
1&0&0&0\\
0&1&0&0\\
0&1&1&0\\
0&0&0&1\\
\end{smallmatrix}\right)\!
\left(\begin{smallmatrix}
1&0&0&0\\
0&1&0&0\\
0&0&-t&0\\
0&0&0&1\\
\end{smallmatrix}\right)
\left(\begin{smallmatrix}
1&0&0&0\\
0&1&0&0\\
0&0&1&-1\\
0&0&0&1\\
\end{smallmatrix}\right)
=
$$
$$
\exp(E_{32})\exp(sE_{33})\exp(-E_{34}),
$$
and
\begin{equation*}
S^2(e^{E_{32}})\!=\!
\left(\begin{smallmatrix}
1&0&0&0&0&0&0&0&0&0\\
0&1&0&0&0&0&0&0&0&0\\
0&1&1&0&0&0&0&0&0&0\\
0&1&0&1&0&0&0&0&0&0\\
0&0&1&0&1&0&0&0&0&0\\
0&0&1&0&2&1&0&0&0&0\\
0&0&0&0&0&0&1&0&0&0\\
0&0&0&0&0&0&0&1&0&0\\
0&0&0&0&0&0&0&1&1&0\\
0&0&0&0&0&0&0&0&0&1\\
\end{smallmatrix}\right),
\,\,\,
S^2(e^{-E_{34}})\!=\!\left(\begin{smallmatrix}
1&0&0&0&0&0&0&0&0&0\\
0&1&0&0&0&0&0&0&0&0\\
0&0&1&0&0&0&0&0&0&0\\
0&0&0&1&0&0&-1&0&0&0\\
0&0&0&0&1&0&0&-1&0&0\\
0&0&0&0&0&1&0&0&-2&1\\
0&0&0&0&0&0&1&0&0&0\\
0&0&0&0&0&0&0&1&0&0\\
0&0&0&0&0&0&0&0&1&-1\\
0&0&0&0&0&0&0&0&0&1\\
\end{smallmatrix}\right)
\end{equation*}
we get
$$
S^2(\sigma_3)\stackrel{\text{Quant}}{\mapsto}$$
$$
\left(\!\begin{smallmatrix}
1&0&0&0&0&0&0&0&0&0\\
0&1&0&0&0&0&0&0&0&0\\
0&1&1&0&0&0&0&0&0&0\\
0&1&0&1&0&0&0&0&0&0\\
0&0&1&0&1&0&0&0&0&0\\
0&0&1&0&(1+q)&1&0&0&0&0\\
0&0&0&0&0&0&1&0&0&0\\
0&0&0&0&0&0&0&1&0&0\\
0&0&0&0&0&0&0&1&1&0\\
0&0&0&0&0&0&0&0&0&1\\
\end{smallmatrix}\!\right)\!\!\!
\left(\!\begin{smallmatrix}
1&0&0&0&0&0&0&0&0&0\\
0&1&0&0&0&0&0&0&0&0\\
0&0&1&0&0&0&0&0&0&0\\
0&0&0&-t&0&0&1&0&0&0\\
0&0&0&0&-t&0&0&0&0&0\\
0&0&0&0&0&qt^2&0&0&0&1\\
0&0&0&0&0&0&1&0&0&0\\
0&0&0&0&0&0&0&1&0&0\\
0&0&0&0&0&0&0&0&-t&0\\
0&0&0&0&0&0&0&0&0&1\\
\end{smallmatrix}\!\right)\!\!\!
\left(\!\begin{smallmatrix}
1&0&0&0&0&0&0&0&0&0\\
0&1&0&0&0&0&0&0&0&0\\
0&0&1&0&0&0&0&0&0&0\\
0&0&0&1&0&0&-1&0&0&0\\
0&0&0&0&1&0&0&-1&0&0\\
0&0&0&0&0&1&0&0&-(1+q)&1\\
0&0&0&0&0&0&1&0&0&0\\
0&0&0&0&0&0&0&1&0&0\\
0&0&0&0&0&0&0&0&1&-1\\
0&0&0&0&0&0&0&0&0&1\\
\end{smallmatrix}\!\right)
$$
$$
=\left(\begin{smallmatrix}
1&0&0&0&0&0&0&0&0&0\\
0&1&0&0&0&0&0&0&0&0\\
0&1&1&0&0&0&0&0&0&0\\
0&0&0&-t&0&0&0&0&0&0\\
0&0&1&0&-t&0&0&0&0&0\\
0&0&1&0&-t(1+q)&qt^2&0&0&0&0\\
0&0&0&0&0&0&1&0&0&0\\
0&0&0&0&0&0&0&1&0&0\\
0&0&0&0&0&0&0&1&-t&0\\
0&0&0&0&0&0&0&0&0&1\\
\end{smallmatrix}\right)
\left(\begin{smallmatrix}
1&0&0&0&0&0&0&0&0&0\\
0&1&0&0&0&0&0&0&0&0\\
0&0&1&0&0&0&0&0&0&0\\
0&0&0&1&0&0&-1&0&0&0\\
0&0&0&0&1&0&0&-1&0&0\\
0&0&0&0&0&1&0&0&-(1+q)&1\\
0&0&0&0&0&0&1&0&0&0\\
0&0&0&0&0&0&0&1&0&0\\
0&0&0&0&0&0&0&0&1&-1\\
0&0&0&0&0&0&0&0&0&1\\
\end{smallmatrix}\right)
$$$$
=\left(\begin{smallmatrix}
1&0&0&0&0&0&0&0&0&0\\
0&1&0&0&0&0&0&0&0&0\\
0&t&1&0&0&0&0&0&0&0\\
0&0&0&-t&0&0&t&0&0&0\\
0&0&1&0&-t&0&0&t&0&0\\
0&0&1&0&-t(1+q)&qt^2&0&t(1+q)&-t^2(1+q)&t^2\\
0&0&0&0&0&0&1&0&0&0\\
0&0&0&0&0&0&0&1&0&0\\
0&0&0&0&0&0&0&1&-t&t\\
0&0&0&0&0&0&0&0&0&1\\
\end{smallmatrix}\right).
$$
Finally,
$$
C_5[S^2(\rho^{(t)}_5(\sigma_3)_e]_qC_5^{-1}=
$$
$$
\left(\begin{smallmatrix}
1&-1&0& 0&0&0   &0 &0 &0&0\\
0& 1&0&-1&0&0   &-1&0 &0&0\\
0&-1&1& 1&-1&0  & 1&-1&0&0\\
0&0 &0& 1&0&0   &-1&0 &0&0\\
0&0 &0&-1& 1&0&1&-1&0 &0\\
0&0 &0& 0&-1&1  &0 &1 &-1&0\\
0&0&0&0&0&0& 1& 0& 0&0\\
0&0&0&0&0&0&-1& 1& 0&0\\
0&0&0&0&0&0&0 &-1& 1&0\\
0&0&0&0&0&0&0 & 0&-1&1\\
\end{smallmatrix}\right)
\left(\begin{smallmatrix}
1&0&0&0&0&0&0&0&0&0\\
0&1&0&0&0&0&0&0&0&0\\
0&t&1&0&0&0&0&0&0&0\\
0&0&0&-t&0&0&t&0&0&0\\
0&0&1&0&-t&0&0&t&0&0\\
0&0&1&0&-t(1+q)&qt^2&0&t(1+q)&-t^2(1+q)&t^2\\
0&0&0&0&0&0&1&0&0&0\\
0&0&0&0&0&0&0&1&0&0\\
0&0&0&0&0&0&0&1&-t&t\\
0&0&0&0&0&0&0&0&0&1\\
\end{smallmatrix}\right)
\times
$$
$$
\left(\begin{smallmatrix}
1&1&0&1 &0&0&1&0&0&0\\
0&1&0&1 &0&0&1&0&0&0\\
0&1&1&1 &1&0&1&1&0&0\\
0&0&0&1 &0&0&1&0&0&0\\
0&0&0&1 &1&0&1&1&0&0\\
0&0&0&1 &1&1&1&1&1&0\\
0&0&0&0&0&0&1&0&0&0\\
0&0&0&0&0&0&1&1&0&0\\
0&0&0&0&0&0&1&1&1&0\\
0&0&0&0&0&0&1&1&1&1\\
\end{smallmatrix}\right)=\left(\begin{smallmatrix}
1&0&0&0 &0&0&0&0&0&0\\
0&0&t&0 &0&0&0&0&0&0\\
0&1&1-t&0 &0&0&0&0&0&0\\
0&0&0  &1 &0&0&0&0&0&0\\
0&0&0&0&0&1  &0&0&0&0\\
0&0&0&0&1&1-t&0&0&0&0\\
0&0&0&0&0&0&1&0&0&0\\
0&0&qt(t-1)&0&0&t(t-1)&0&qt^2&0&t(t-1)\\
0&0&0&0&0&0&0&0&1&t\\
0&0&0&0&0&0&0&0&0&1-t\\
\end{smallmatrix}\right)
\stackrel{(\ref{Kram-Big-new})}{=}k_5^{(t,q)}(\sigma_3).
$$
Similarly, we can prove (\ref{K-eq-B(n)}) for general $B_n$.
%\qed
%\end{pf}%

%%%%%%%%%%%%%%%%%%%%%%%
%\newpage

\section{How to construct new representations of the braid group $B_n$}
\subsection{Quantization of the $m^{th}$ symmetric power of the Burau representation}
We start with the well-known representation (in fact, reduced Burau  representation for $t=-1$)
$T_{32}:=\rho^{(-1)}_3:B_3\rightarrow {\rm SL}(2,{\mathbb Z})$
defined by (\ref{Bur(3)})
\begin{equation}
\label{Bur(3)1}
 \sigma_1\mapsto \left(\begin{smallmatrix}
1&-1\\
0&1\\
\end{smallmatrix}\right)
\quad
\sigma_2\mapsto  \left(\begin{smallmatrix}
1&0\\
1&1\\
\end{smallmatrix}\right),\quad\text{or}\quad
 \sigma_1\mapsto \left(\begin{smallmatrix}
1&1\\
0&1\\
\end{smallmatrix}\right)
\quad
\sigma_2\mapsto  \left(\begin{smallmatrix}
1&0\\
-1&1\\
\end{smallmatrix}\right).
\end{equation}
This representation is just $\exp$ of the fundamental (or natural) representations $\pi_2$ of
the Lie algebra $\mathfrak{sl}_2$ (see (\ref{B_3-1}) and (\ref{B_3-2})).
\begin{rem}
\label{r.Rep(B_n)-general}
We have {\it several possibility} to generalize representation $T_{32}$ of the group $B_3$.
{\it First}, we can generalize $T_{32}$ for an arbitrary $B_n$ (denote by $T_{n2}$), {\it second}, we can generalize $T_{32}$
of the same group $B_3$ but for an arbitrary dimension $m$ (denote by $T_{3m}$), and {\it finally}, we can generalize representations $T_{n2}$ of $B_n$ for higher dimension $m$, (denote by $T_{nm}'$) and extend representation $T_{3m}$ of $B_3$ for other groups $B_n$  (denote by $T_{nm}''$).
Finally, we  have two family of representations $$T_{nm}',T_{nm}'':B_n\rightarrow {\rm GL}(m,\mathbb C).$$
\end{rem}
%
%\begin{rem}
%In two words, $T_{n2}$ is the Burau representation of $B_n$, and $T_{3n}''$ is the Lawrence-Krammer representation of $B_n$.
%\end{rem}
%
To be more precise, the {\it first step} gives us all Burau's representations
$\rho^{(t)}_{n+1}:B_{n+1}\to {\rm GL}_{n}({\mathbb Z}[t,t^{-1}])$, due to Theorem~\ref{l.Bur-sl_n},
if we start with $\rho^{(t)}_{3}$ and  take for $\pi$ the natural representation $\pi_n$ of $\mathfrak{sl}_n$.

On the {\it second step}, we obtain the Humphry representation of $B_3$ in
all dimension $m+1$ as the $m^{th}$ symmetric power of the Burau representation
$\rho_3^{-1}$, due to Theorem~\ref{B_3-sl_2}. Remarkably, that {\it quantization} and {\it diagonal deformation} (governed by the matrix $\Lambda_n={\rm diag}(\lambda_r)_{r=0}^n$, see Remark~\ref{pi^Lambda_3})
of the Humphry representations  (see (\ref{Rep(q)})) {\it include all irreducible representations} of $B_3$ in
dimension $\leq 5$, due to results of Tuba and Wensel \cite{Tub01,TubWen01} and our equivalent description of their representations 
\cite{KosAlb07q,Kos07q}, see Theorem~\ref{Rep(B_3)<5}. Our approach using $q$-Pascal triangle allows us, for free, extend representation $T_{32}$ of $B_3$ for arbitrary dimension (see Theorem~\ref{t.Rep(q)}).

On the {\it third step}, we obtain first, the Lawrence-Krammer representation as quantization of the symmetric square of the Burau represention,
Theorem~\ref{Kr=S^2_q(Bur)}. Due to formulas (\ref{T(33)1}), (\ref{T(33)2}) and Theorem~\ref{Kr=S^2_q(Bur)}, the Lawrence-Krammer representations
of $B_3$ describe all irreducible representations, up to the scalar factor, in dimensions $3$.
\begin{rem}
\label{r.(Symm^3)_q}
%So, to write down new representations for all $B_n$ we should start from representations of $B_3$ in dimension $n\geq 4$. In this case we have
It would be nice to give explicit formulas  for quantization of  the symmetric cube and  higher symmetric power of the Burau representation.
\end{rem}
\begin{table}[ht]
  \centering
  %\caption{Caption for the table.}
  %\label{tab:table1}
  \begin{tabular}{|c|c|c|c|}
  \hline
  $V\!={\mathbb C}^n$&$B_3$&  $B_n$         &type of representation         \\
  \hline
  $2$&$T_{3,2}^{\Lambda_1,1,\pi_2}$&$\rho_n^{(t)}$&reduced Burau representation\\
  %&$T_{n,2}^{\pi_{n-1}}$ \\
  \hline
  $3$&$T_{3,3}^{\Lambda_2,q,\pi_2}$&$k^{(t,q)}_n$&Lawrence--Krammer representation\\
  %&$T_{n,3}^{\Lambda,q,\pi_{n-1}}$ \\
  \hline
  $ 4$&$T_{3,4}^{\Lambda_3,q,\pi_2}$&$[S^3(\rho_{n}^{(t)})]_q$ &Symmertic cube\\
  %& $T_{n,4}^{\Lambda,q,\pi_{n-1}}$ \\
    \hline
  $ 5$&$T_{3,5}^{\Lambda_4,q,\pi_2}$& $[S^4(\rho_{n}^{(t)})]_q$         &$4^{th}$ Symmertic power\\
  %& $T_{n,5}^{q,\pi_{n-1}}$  \\
    \hline
   $m+1$&$T_{3,m+1}^{\Lambda_m,q,\pi_2}$& $[S^m(\rho_{n}^{(t)})]_q$         &$m^{th}$ Symmertic power\\
   %&$T_{n,m+1}^{q,\pi_{n-1}}$ \\
  \hline
  \end{tabular}
\end{table}
\subsection{Representations of a Lie algebra $\mathfrak{sl}_{n+1}$ }
We recall, following to J.P.Serre \cite{Serr64}, Part I, Ch. VII, \S 1; Part III, Ch. VII, \S 6,
that all  irreducible finite-dimensional representations of a {\it simple} Lie algebra $\mathfrak{sl}_{n+1}$  are {\it highest weight representations} and are contained in tensor powers of the {\it natural representation} $\pi_{n+1}$ \cite{Weyl46}.

The {\it Cartan subalgebra} $\mathfrak h$ of $\mathfrak{sl}_{n+1}$  consists of all matrices of the form $H={\rm diag}(\lambda_1,\dots, \lambda_{n+1})$ with $\sum_k\lambda_k=0$. The {\it roots} of a Lie algebra $\mathfrak{sl}_{n+1}$ are linear functionals $\alpha_{i,j}\in {\mathfrak h}^*,\,\,i\not=j,$ of the form
$$
\alpha_{i,j}(H)=\lambda_i-\lambda_j,\quad \text{where}\quad  H={\rm diag}(\lambda_1,\dots, \lambda_{n+1}).
$$
The {\it fundamental wights} $\omega_k$ are defined by
$$
\omega_k(H)=\lambda_1+\dots+\lambda_k,\quad \text{where}\quad  H={\rm diag}(\lambda_1,\dots, \lambda_{n+1}).
$$
The fundamental weight $\omega_1$ is the fundamental weight of the {\it standard representation} $\pi_{n+1}$ of the Lie algebra $\mathfrak{sl}_{n+1}$ in $E\!=\!{\mathbb C}^{n+1}$. Moreover, the  fundamental weight $\omega_k$ is the highest weight of the representation
$\wedge^k(\pi_{n+1})$ of $\mathfrak{sl}_{n+1}$ in the space $\wedge^kE$ (see \cite[1.7.4.6]{Serr64}), where $\wedge^k(A)$ is $k^{th}$ {\it exteriour power} of $A$.

In particular case, all irreducible representations of $\mathfrak{sl}_2$ are ${\rm Sym}^n(\pi_2)$ in ${\rm Sym}^n({\mathbb C}^2),\,\,n\in {\mathbb N}$. For $\mathfrak{sl}_3$  all irreducible representations  are contained in  tensor product of $\pi_3$ and $\wedge^2(\pi_3)$, corresponding t fundamental weights $\omega_1$ and $\omega_2$.
All irreducible representations of $\mathfrak{sl}_{n+1}$ are contained in  tensor products of  $\wedge^k(\pi_{n+1}),\,\,1\leq k\leq n$.

Let $\pi_{n}\!:\!\mathfrak{sl}_n\rightarrow {\rm End}({\mathbb C}^m)$ be any finite-dimensional representation of a Lie algebra $\mathfrak{sl}_n$ then
$\exp(\pi_{n}):{\rm SL}(n,\mathbb C)\rightarrow {\rm GL}(m,\mathbb C)$ is representation of ${\rm SL}(n,\mathbb C)$.
We can define representations of $B_{n+1}$ as follows
\begin{equation}
\label{Rep(B_n,pi_n)}
\sigma_k\to S^m\big(\exp(\pi_{n})\circ \rho_n^{(t)}(\sigma_k)\big),\quad
\sigma_k\to \wedge^m\big(\exp(\pi_{n})\circ \rho_n^{(t)}(\sigma_k)\big).
\end{equation}
Consider the exteriour square $\wedge^2(\rho_{2,4}^{(t)})$ of the Burau representation $\rho_{2,4}^{(t)}$ of the group
$B_4$ defined by (\ref{Bur_n-conj}) and (\ref{Bur_n(b)-conj}) in the basis $e_{ij}^\wedge,\,\,1\leq i<j\leq 3$ of the space $\wedge^2({\mathbb C}^3)$
$$
e_{12}^{\wedge}=e_1\otimes e_2-e_2\otimes e_1,\quad
e_{13}^{\wedge}=e_1\otimes e_3-e_3\otimes e_1,\quad
e_{23}^{\wedge}=e_2\otimes e_3-e_3\otimes e_2.
$$
\begin{lem}
\label{ext^2-Bur-sim-Bur}
Representation $\wedge^2(\rho_{2,4}^{(t)})$ of the group $B_4$
is defined by
%the following formulas :
\begin{equation}
\label{ext(2)-Bur_4}
\sigma_1\to\left(\begin{smallmatrix}
-t&0&0\\
0&-t&t\\
0&0&1
\end{smallmatrix}\right),\quad
\sigma_2\to\left(\begin{smallmatrix}
-t&t&0\\
0&1&0\\
0&1&-t
\end{smallmatrix}\right),\quad
\sigma_3\to\left(\begin{smallmatrix}
1&0&0\\
1&-t&0\\
0&0&-t
\end{smallmatrix}\right).
\end{equation}We have
\begin{equation}
\label{Ext-Bur-Bur}
\wedge^2\rho_{2,4}^{(t)}\not\sim \rho_4^{(t)},\quad\text{for}\,\,\, t\not=-1,\quad
\text{but}\,\,\, \wedge^2\rho_{2,4}^{(t)}=-t S_3\rho_{1,4}^{(-t^{-1})}S_3
\end{equation}
where $\rho_{1,4}^{(t)}$ is  defined by (\ref{r-Bur_n}) and (\ref
{r-Bur_n(b)}) and  $S_3=\sum_{k=1}^2E_{k,4-k}=
\left(\begin{smallmatrix}
0&0&1\\
0&1&0\\
1&0&0
\end{smallmatrix}\right).$
\end{lem}
%
%\begin{pf} 
%
The first part of (\ref{Ext-Bur-Bur})
follows from comparison of two spectra: ${\rm Sp}\rho_{24}^{(t)}(\sigma_1)=\{-t,1,1\}$ and
${\rm Sp}\big(\wedge^2\rho_{24}^{(t)}(\sigma_1)\big)=\{-t,1,-t\}=-t\{1,-t^{-1},1\}.$ The second part follows from
(\ref{r-Bur_n}) and (\ref{r-Bur_n(b)}).
%\qed
%\end{pf}

%
\begin{prb}
\label{r.Rep(B_n,pi_n)_q}
To quantize the representation of $B_n$ defined by (\ref{Rep(B_n,pi_n)}).
\begin{equation}
\label{Rep(B_n,pi_n)_q}
\sigma_k\to \big[S^m\big(\exp(\pi_{n})\circ \rho_n^{(t)}(\sigma_k)\big) \big]_q,\quad
\sigma_k\to \big[\wedge^m\big(\exp(\pi_{n})\circ \rho_n^{(t)}(\sigma_k)\big) \big]_q.
\end{equation}
\end{prb}

\section{Some polynomial invariants for knots}
\label{s.3}
\subsection{Some polynomial invariants for the trefoil, see
{\rm \cite{BarMor}}}
%,Wei-Wolfram

%  and  {\rm \cite{Wei-Wolfram}}}
%%%%%%%%%%%
\iffalse
{\rm \cite{BarMor}}
The {\it Alexander polynomial} of the {\it trefoil knot} is
\begin{equation}
\Delta(t)=t-1+t^{-1}=t^{-1}(t^2-t+1),
\end{equation}
and the {\it Conway polynomial} is
\begin{equation}
\nabla(z)=z^2+1,
\end{equation}
The {\it Jones polynomial} is
\begin{equation}
V(q)=q^{-1}+q^{-3}-q^{-4}=q^{-4}(q^3+q-1),
\end{equation}
and the {\it Kauffman polynomial} of the trefoil is
\begin{equation}
L(a,z)=za^5+z^2a^4-a^4+za^3+z^2a^2-2a^2.
\end{equation}
\fi
%%%%%%%%%%%%%%
The Alexander polynomial $\Delta(x)$, BLM/Ho polynomial $Q(x)$, Conway polynomial $\nabla(x)$, HOMFLY polynomial $P(l,m)$, Jones polynomial $V(t)$, and Kauffman polynomial $F(a,z)$ of the trefoil knot are
\begin{eqnarray}
\Delta(x)	=	x-1+x^{-1},\\	
Q(x)	=	2x^2+2x-3,\\
\nabla(x)	=	x^2+1,\\		
P(l,m)	=	-l^4+m^2l^2-2l^2,\\		
V(t)_{3_1}	=	t+t^3-t^4,\,\,V(t)_{3_2}	=	t^{-1}+t^{-3}-t^{-4},\\	
F(a,z)	=	-a^4-2a^2+(a^4+a^2)z^2+(a^5+a^3)z.	
\end{eqnarray}

\subsection{The Burau representation and  the Alexander polynomial}

Let $K$ be some link, by {\it Alexander's theorem} this link can be
obtained as the closure of some braid $X\in B_{n}$, notation $K={\it
Cl}(X)$. The {\it Alexander polynomial} $\Delta_K(t)$ can be
obtained using the reduced Burau representation by the following
formula, \cite[(4.3), p.64]{BirBre93}:
\begin{equation}
\label{Burau-alex_pol} \Delta_{{\it Cl}(X)}(t)\!=\!\frac{(1-t){\rm
det}(\rho_n^{(t)}(X)-I_{n-1})}{1-t^n}
%{1+t+...+t^{n-1}}
\!=\!
\frac{{\rm
det}(\rho_n^{(t)}(X)-I_{n-1})}{{\rm
det}(\rho_n^{(t)}(\sigma_1\dots \sigma_{n-1})-I_{n-1})}.
\end{equation}

\subsection{The Krammer representation and  the corresponding rational function}
\begin{df}
\label{d.Krammer-pol}
By analogy with the definition of the Alexander polynomial we can introduce the following
%polynomial
rational function $k_n(t,q)$  connected with the Lawrence-Krammer representation:
\begin{equation}
\label{Kram-pol}
 k_{n,{\it Cl}(X)}(t,q)=
\frac{{\rm det}(k_n^{(t,q)}(X)-I_{n(n-1)/2})}
{{\rm
det}(k_n^{(t,q)}(\sigma_1\dots \sigma_{n-1})-I_{n(n-1)/2})}.
\end{equation}
\end{df}
\begin{rem} Recall  that some polynomial $P(x,y,\dots)$ defined on the set of all  knots is {\it knot invariant} if and only if it respect the
the {\it first and the second Markov moves} (see, e.g., \cite[2.5. Markov's theorem, p.67]{KasTur08}), i.e., for any $X\in B_n$ the following properties hold:
\begin{equation}
\label{mark1}
P_{{\it Cl}(wXw^{-1})}(x,y,\dots)=P_{{\it Cl}(X)}(x,y,\dots)\quad\text{for all}\quad w\in B_n,
\end{equation}
\begin{equation}
\label{mark2}
P_{{\it Cl}(X\sigma_{n}^{\pm 1})}(x,y,\dots)=P
_{{\it Cl}(X)}(x,y,\dots).
\end{equation}
\end{rem}
\begin{thm}
\label{t.Kram-pol}
The
%polynomial
rational function $k_{n,{\it Cl}(X)}(t,q)$ defined by (\ref{Kram-pol}) respect the first Markov's move, i.e.,
depends only on conjugacy classes in $B_n$.
%is the knot invariant.
\end{thm}
%
%\begin{pf}
%
We show that
\begin{equation}
\label{Mark1}
k_{{\it Cl}(wXw^{-1})}(t,q)=k_{{\it Cl}(X)}(t,q)\quad\text{for all}\quad w\in B_n,
\end{equation}
Since ${\rm det}(AB)={\rm det}(A){\rm det}(B)$ we get by definition
$$
{\rm det}\big(k_n^{(t,q)}( wXw^{-1})-I_{n(n-1)/2}\big)=
{\rm det}\big(k_n^{(t,q)}( wXw^{-1})-k_n^{(t,q)}( ww^{-1})\big)
$$
$$
={\rm det}\Big(k_n^{(t,q)}( w)\big( k_n^{(t,q)}(X)-I_{n(n-1)/2}\big)
k_n^{(t,q)}( w)^{-1}\Big)
$$
$$
={\rm det}(k_n^{(t,q)}(X)-I_{n(n-1)/2}),
$$
that proves (\ref{Mark1}).
%\qed
%\end{pf}
%
\begin{prb}
Study the  behaviour   of the rational function
%polynomial
$k_n(t,q)$ with respect to the second Markov's move.
\end{prb}
\subsection{Alexander polynomial for the trefoil}
To calcultate the Alexander polynomial for the trefoil knot $3_1\!=\!{\it Cl}(\sigma_1^3)\!=\!{\it Cl}(\sigma_1^3\sigma_2)$
we get for the representation $\rho_2^{(t)}(\sigma_1)=-t$ of the group $B_2$
\begin{equation}
\label{Alex1}
\Delta_{3_1}(t)=\frac{{\rm det}\big(\rho_2^{(t)}(\sigma_1^3)-I_1\big)}{{\rm det}\big(\rho_2^{(t)}(\sigma_1)-I_1\big)}=\frac{-t^3-1}{-(t+1)}=t^2-t+1.
\end{equation}
{\it The second Markov move (see (\ref{mark2}))}. For the Burau representation $\rho_3^{(t)}$ of the group  $B_3$
$$
\rho_3^{(t)}(\sigma_1)=\left(\begin{smallmatrix}
-t&t\\
0&1\\
\end{smallmatrix}\right),\quad
\rho_3^{(t)}(\sigma_2)= \left(\begin{smallmatrix}
1&0\\
1&-t\\
\end{smallmatrix}\right)
$$
we have
$$
\rho_3^{(t)}(\sigma_1\sigma_2^3)=\left(\begin{smallmatrix}
-t&t\\
0&1\\
\end{smallmatrix}\right) \left(\begin{smallmatrix}
1&0\\
1&-t\\
\end{smallmatrix}\right)^3=
 \left(\begin{smallmatrix}
t^3-t^2&-t^4\\
t^2-t+1&-t^3\\
\end{smallmatrix}\right).
$$
Hence,
$$
{\rm
det}\big(\rho_3^{(t)}(\sigma_1\sigma_2^3)-I_2\big)=
{\rm
det}\left(\begin{smallmatrix}
t^3-t^2-1&-t^4\\
t^2-t+1&-t^3-1\\
\end{smallmatrix}\right)=
$$
$$
-(t^3+1)(t^3-t^2-1)+t^4(t^2-t+1)=t^4+t^2+1,
$$
and we obtaine the same result as in (\ref{Alex1}):
\begin{equation}
\label{Alex2}
\Delta_{3_1}(t)=
\frac{{\rm det}\big(\rho_3^{(t)}(\sigma_1^3\sigma_2)-I_2\big)}{{\rm det}\big(\rho_3^{(t)}(\sigma_1\sigma_2)-I_2\big)}
=\frac{t^4+t^2+1}{t^2+t+1}=t^2-t+1.
\end{equation}

\subsection{Krammer representation and the trefoil}
The Krammer representation for the group $B_2$ is as follows (see (\ref{Kram-Big-new}))
%\begin{equation}
$\sigma_1\to t^2q,$
%\end{equation}
therefore, using (\ref{Kram-pol}) we get for the trefoil knot $3_1\!=\!{\it Cl}(\sigma_1^3)$
(see (\ref{Alex1}))
\begin{equation}
\label{Kram-3(1)}
k_{2,3_1}(t,q)=\Delta_{3_1}(t^2q)=(t^2q)^2-t^2q+1=t^4q^2-t^2q+1.
\end{equation}
\begin{rem}
We have
$$
k_{2,3_1}(1,q)=q^2-q+1=\Delta_{3_1}(q),\quad k_{2,3_1}(t,1)=t^4-t^2+1=\Delta_{3_1}(t^2).
$$
\end{rem}
Calculate now $k_{3,3_1}(t,q)$ for the trefoil knot $3_1={\it Cl}(\sigma_1^3\sigma_2)$ using formulas (\ref{T(33)(1,2)}),
(\ref{S^2(s1)v]q}), (\ref{S^2(s2)v]q}) and (\ref{kra-B3}).
\begin{lem}
We have
\begin{equation}
\label{det1}
{\rm det}\big(k^{(t,q)}_3(\sigma_1\sigma_2)-I_3\big)=
%%%%%%%%%%%%%%%
\iffalse
\left|\begin{smallmatrix}
-1&0&t^4q\\
0&-tq^2-1&t^3q\\
1&-t(1+q)& t^2q-1\\
\end{smallmatrix}\right|=
\fi
%
t^6q^2-1,
\end{equation}
%
%Further, using (\ref{s_1^3s_2}) we obtain
\begin{equation}
\label{det2}
{\rm det}\big(k^{(t,q)}_3(\sigma_1^3\sigma_2)-I_3\big)=
%
%t^{12}q^5-t^6q^3+t^6q^2-1=(t^6q^3+1)(t^6q^{2}-1).
t^{12}q^4-t^6q(1-q)(1-t)(1-tq)-1,
\end{equation}
%Finally, we get
\begin{equation}
\label{det1/2}
k_{3,3_1}(t,q)=\frac{t^{12}q^4-t^6q(1-q)(1-t)(1-tq)-1}{ t^6q^2-1}.
\end{equation}
\end{lem}
\begin{rem}
\label{2Marc-move}
The rational function
%polynomial
$k_n(t,q)$ does not respect the second Markov's move.
% but gives some factors.
Indeed, in the particular case of $q=1$ we get
\begin{equation}
k_{3,3_1}(t,1)=\frac{t^{12}-1}{t^6-1}=t^6+1=(t^2+1)(t^4-t^2+1)= (t^2+1)  k_{2,3_1}(t,1).
\end{equation}
\end{rem}
\begin{rem}
 We have for $q=1$
\begin{equation}
\label{}
k_{3,3_1}(t,1)= (t^2+1)  k_{2,3_1}(t,1)=(t^2+1)\Delta(t^2),\quad \text{\it the Alexander polynomial}.
\end{equation}
In the particular case of $t=1$ we get
\begin{equation}
k_{3,3_1}(1,q)=\frac{q^4-1}{q^2-1}=q^2+1=\nabla(q),\quad \text{\it the Conwey polynomial}.
\end{equation}
\end{rem}
%
%\begin{pf}
%
To prove (\ref{det1}) and  (\ref{det2}), using (\ref{T(33)(1,2)}) we have
\begin{equation}
\label{T(33)(1,2)1}
%[S^2(\sigma_1)_e]_q:=
\sigma_1\mapsto
\left(\begin{smallmatrix}
t^2q&-t^2(1+q)&t^2\\
0&-t&t\\
0&0&1
\end{smallmatrix}\right)
,\quad
%[S^2(\sigma_2)_e]_q:
\sigma_2\mapsto
\left(\begin{smallmatrix}
1&0&0\\
1&-t&0\\
1&-t(1+q)&t^2q
\end{smallmatrix}\right),
\end{equation}
therefore
$$
\sigma_1\sigma_2
\mapsto\left(\begin{smallmatrix}
0&0&t^4q\\
0&-tq^2&t^3q\\
1&-t(1+q)& t^2q\\
\end{smallmatrix}\right):=B(t,q).
$$
To calculate ${\rm det}(B(t,q)-I_3)$ we use the following lemma.
\begin{lem}
[Lemma 1.4.5 in \cite{Kos_B_09}]
\label{l.g-ch-pol}
%\label{l.detC-LI}
For the {\rm generalized characteristic po\-lynomial} $P^g_C(\lambda)$ of a matrix
$C\!\in\!{\rm Mat}(m,{\mathbb C})$ defined by (\ref{g-ch-pol}) and
$\lambda
%=(\lambda_k)_{k=1}^m 
\in {\mathbb C}^m$ we have 
\begin{equation}
\label{g-ch-pol}
P^g_C(\lambda)= {\rm
det}\,\Big(C+\sum_{k=1}^m\lambda_kE_{kk}\Big)\!=
\end{equation}
\begin{equation}
\label{detC-LI}
\!{\rm det}\,C+
\sum_{r=1}^m\sum_{1\leq i_1<i_2<...<i_r\leq
m}\lambda_{i_1}\lambda_{i_2}...\lambda_{i_r}A^{i_1i_2...i_r}_{i_1i_2...i_r}(C),
\end{equation}
where $A^{i_1i_2...i_r}_{i_1i_2...i_r}(C)$ are cofactors of $C$.
\end{lem}
%%
\iffalse
As the particular case, we obtain the well-known formula for the characteristic polynomial $P_C(t)={\rm det}(tI-C)$ of the
matrix $C\!\in\!{\rm Mat}(m,{\mathbb C})$
\begin{equation}
\label{P_C(t)}
P_C(t)=\sum_{k=0}^{m}t^{m-k}(-1)^k{\rm tr}(\Lambda^kC)=
\end{equation}
\begin{equation}
\label{P_C(t)1}
-\Big({\rm det}\,C+
\sum_{r=1}^m(-t)^r\sum_{1\leq i_1<i_2<...<i_r\leq
m}A^{i_1i_2...i_r}_{i_1i_2...i_r}(C)\Big),
\end{equation}
where ${\rm tr}(\Lambda^kC)$ is the trace of the $k^{th}$ exterior power of $C$ which has
 dimension $\left(\begin{smallmatrix}n\\2\end{smallmatrix}\right)$. This trace may be computed as the sum of all principal minors of $C$ of size $k$.
%
\fi
%%%%
By Lemma~\ref{l.g-ch-pol} we have for $A=B(t,q)$ and $\lambda_k=-1,\,\,k=1,2,3$
$$
{\rm det}(B(t,q)-1)={\rm det}A(t)-(A^1_1+A^2_2+A^3_3)+(A^{12}_{12}+A^{13}_{13}+A^{23}_{23})-A^{123}_{123}=
$$
$$
t^6q^2-(t^4q-t^4q)+(t^2q-t^2q)-1=t^6q^2-1.
$$
To prove (\ref{det2}) we have
$$
\sigma^3_1\sigma_2=\sigma^2_1\sigma_1\sigma_2\mapsto
\left(\begin{smallmatrix}
t^4q^2&t^3(1+q)(1-tq)&t^2[t^2q-t(1+q)+1]\\
0&t^2&-t(t-1)\\
0&0& 1\\
\end{smallmatrix}\right)
\left(\begin{smallmatrix}
0&0&t^4q\\
0&-tq^2&t^3q\\
1&-t(1+q)& t^2q\\
\end{smallmatrix}\right)=
$$
\begin{equation}
\label{s_1^3s_2}
\left(\!\begin{smallmatrix}
t^2[t^2q-t(1+q)+1]&(1+q)t^3(tq-1)(t^2q-t+1)&t^4q[t^4q^2+t^2(1+q)(1-tq)+t^2q-t(1+q)+1]\\
-t(t-1)&-t^4q+t^2(1+q)(t-1)&t^3q(t^2-t+1)\\
1&-t(1+q)& t^2q\\
\end{smallmatrix}\!\right)=\!:A(t,q).
\end{equation}
For $A=A(t,q)$ and $\lambda_k=-1,\,\,k=1,2,3$ we get
$$
{\rm det}(A(t,q)-1)={\rm det}A(t,q)-(A^1_1+A^2_2+A^3_3)+(A^{12}_{12}+A^{13}_{13}+A^{23}_{23})-A^{123}_{123}
$$
$$
=(t^3q)^4-t^6q\big(
1-[t^2q^2-q(1+q)t+(1+q)]+[t^2q-(1+q)t+1]
\big)
$$
$$
+t^2\big(
[t^2q-(1+q)t+1]-[t^2q-(1+q)t+(1+q)]+q
\big)
-1=
$$
$$
t^{12}q^4-t^6q(1-q)[t^2q-(1+q)t+1]-1=t^{12}q^4-t^6q(1-q)(1-t)(1-tq)-1.
$$
%\qed
%\end{pf}

{\bf Acknowledgements 2017.} The author expresses his deep gratitude to
the Max--Planck--Institute of Mathematics for the financial grant and the hospitality in 2016--2017.
Sergei Chmutov have read the first version of this article and made a lot of useful remarks, that improved considerable the content. I am very grateful to him for his valuable help.

\vskip 0.3 cm
{\bf Acknowledgements 2023.}
The author is very grateful to Prof. K.-H. Neeb, Prof. M.~Smirnov and 
Dr P.~Moree
for their personal efforts  to make academic  stays possible at their respective institutes. %A.~Kosyak
The author
visited: MPIM from March to April 2022 and from January to April 2023,
University of Augsburg from June to July 2022, and  University of Erlangen-Nuremberg 
from August to December 2022, all during the Russian invasion in Ukraine.
Also, Prof. R.~Kashaev  kindly invited 
%A.~Kosyak 
him to Geneva.

Further, he 
would like  to pay his respect to Prof. P. Teichner at MPIM, for his
immediate efforts started to help mathematicians
in Ukraine  after the Russian invasion.

Since the spring of 2023 A.~Kosyak is an Arnold Fellow at the
London Institute for Mathematical Sciences, and he would like to express
his gratitude to  Mrs S.~Myers Cornaby   to Miss A.~Ker Mercer and to Dr M.~Hall 
%Madeleine 
and especially to the Director of LIMS Dr T.~Fink and Prof. Y.-H.~He.
%
%\newpage

\end{document}